\documentclass[11pt]{article}

\usepackage[T1]{fontenc}
\usepackage[utf8]{inputenc}
\usepackage{lmodern}

\usepackage{amsmath,amssymb,amsthm,bm}

\usepackage[a4paper,margin=1.15in]{geometry}
\usepackage{graphicx}
\usepackage{caption}
\usepackage{booktabs,array}

\usepackage[colorlinks=true,
linkcolor=blue,
citecolor=blue,
urlcolor=blue]{hyperref}

\usepackage{titlesec}
\titleformat{\section}
{\normalfont\Large\bfseries}{\thesection.}{0.75em}{}
\titleformat{\subsection}
{\normalfont\large\bfseries}{\thesubsection.}{0.75em}{}
\titleformat{\subsubsection}
{\normalfont\normalsize\bfseries}{\thesubsubsection.}{0.75em}{}

\usepackage{abstract}

\newtheorem{theorem}{Theorem}[section]

\newtheorem{definition}[theorem]{Definition}
\newtheorem{remark}[theorem]{Remark}

\newtheorem{condition}[theorem]{Condition}
\newtheorem{corollary}[theorem]{Corollary}

\newcommand{\Sd}{\mathbb{S}^d}

\title{\vspace{-1em}
	Spectral Bayesian Regression on the Sphere
	\vspace{-0.5em}}

\author{
	Claudio Durastanti\\
	\small Department of Basic and Applied Sciences for Engineering\\
	\small Sapienza University of Rome\\
	\small\texttt{claudio.durastanti@uniroma1.it}
}

\date{\small \today}
\begin{document}
	\maketitle
	
	\begin{abstract}
We develop a fully intrinsic Bayesian framework for nonparametric regression on
the unit sphere based on isotropic Gaussian field priors and the harmonic
structure induced by the Laplace--Beltrami operator.
Under uniform random design, the regression model admits an exact diagonalization
in the spherical harmonic basis, yielding a Gaussian sequence representation with
frequency-dependent multiplicities.

Exploiting this structure, we derive closed-form posterior distributions,
optimal spectral truncation schemes, and sharp posterior contraction rates under
integrated squared loss.
For Gaussian priors with polynomially decaying angular power spectra, including
spherical Mat\'ern priors, we establish posterior contraction rates over Sobolev
classes, which are minimax-optimal under correct prior calibration.

We further show that the posterior mean admits an exact variational
characterization as a geometrically intrinsic penalized least-squares estimator,
equivalent to a Laplace--Beltrami smoothing spline.
\vspace{0.5em}

	\textbf{Keywords:}
Bayesian nonparametric regression;
Gaussian process regression;
Spherical harmonics;
Laplace--Beltrami operator;
Posterior contraction;
Mat\'ern priors;
Spectral regularization;
Smoothing splines on manifolds.

\vspace{0.5em}
\textbf{MSC (2020):}
62G08; 62G20; 60G15; 58J35; 42C10

\end{abstract}

\section{Introduction}

Gaussian process (GP) regression is a fundamental tool in nonparametric
statistics, providing a probabilistic framework for function estimation and
uncertainty quantification (see, e.g., \cite{rasmussen2006,wahba90}).
While the classical theory is largely developed for Euclidean domains, many
statistical problems involve data indexed on non-Euclidean spaces, most notably
the unit sphere $\mathbb{S}^d$.
Such settings arise naturally in global geostatistics, climate science,
cosmology, and directional data analysis, where respecting the intrinsic
geometry of the domain is essential for both modeling and inference (see, for
instance, \cite{Gneiting13,MP11}).
Sharp posterior contraction rates for Gaussian process priors in nonparametric
regression and inverse problems have been extensively studied in Euclidean
settings; see, for example,
\cite{CastilloNickl14,vanDerVaartvanZanten08}.

Gaussian processes on $\mathbb{S}^d$ are most naturally constructed using the
intrinsic geometric and harmonic structure of the sphere.
In particular, isotropic Gaussian random fields on $\mathbb{S}^d$ admit spectral
representations in terms of spherical harmonics, which are the eigenfunctions
of the Laplace--Beltrami operator.
Their covariance structure is fully characterized by an angular power spectrum
\cite{MP11,Yadrenko83}.
This spectral viewpoint provides a principled alternative to Euclidean kernel
constructions and underlies both the classical theory of spherical random fields
and modern statistical models for global data \cite{jeong17}.
Moreover, the harmonic characterization enables a precise analysis of sample
regularity, efficient simulation schemes, and intrinsic stochastic partial
differential equation representations of spherical Gaussian fields
\cite{LangSchwab15}.

Gaussian random fields on $\mathbb{S}^d$ provide the probabilistic foundation for
Gaussian process regression on spherical domains, which is obtained by
conditioning such fields on (possibly noisy) observations.
In Euclidean settings, it is well known that the posterior mean of a Gaussian
process coincides with the solution of a penalized least-squares problem in a
reproducing kernel Hilbert space, and, conversely, that a broad class of
classical smoothers can be interpreted as posterior means under suitable
Gaussian priors \cite{kimeldorf70,wahba90}.
This equivalence between Gaussian-process priors and quadratic regularization
extends beyond Euclidean domains and applies to geometric settings and compact
manifolds, where intrinsic differential operators and harmonic structures play a
central role \cite{lindgren2011}.

Our analysis shows that isotropic Gaussian field priors on $\mathbb{S}^d$ induce
an exact Gaussian sequence representation, with harmonic multiplicities growing
as $\ell^{d-1}$, a geometric feature that fundamentally alters the bias–variance
trade-off relative to Euclidean domains.
This viewpoint clarifies connections between Bayesian regression, variational
regularization, kernel methods on manifolds, and Whittle-type spectral inference,
and shows that sharp theoretical results on $\mathbb{S}^d$ require intrinsically
geometric arguments.

While Gaussian process regression, Mat\'ern priors, and spherical harmonics have
all been studied previously, a unified Bayesian theory of nonparametric
regression on the sphere that fully exploits the intrinsic harmonic structure
and yields sharp frequentist guarantees has been lacking.
Existing approaches have either focused on likelihood-based inference for
Gaussian random fields, or on adaptive multiscale methods, where posterior
distributions and exact variational characterizations are typically not
available.

The unit sphere $\mathbb{S}^d$ is a canonical non-Euclidean domain in statistics,
with strong mathematical and applied relevance in geostatistics, climate science,
cosmology, and directional data analysis \cite{Gneiting13,MP11}.
From a theoretical perspective, it is one of the very few compact manifolds that
simultaneously admit a transitive symmetry group, a natural notion of isotropy,
and a fully explicit harmonic decomposition.
As a consequence, isotropic Gaussian fields are exactly diagonal in the spherical
harmonic basis, yielding an exact finite-sample Gaussian sequence representation
under uniform random design.
On general manifolds, where eigenfunctions are not explicit and multiplicities
are irregular, such exact diagonalization and closed-form posterior analysis are
typically unavailable.

The present work provides a fully intrinsic Bayesian framework for nonparametric
regression on the unit sphere under random design.
The assumption of uniform random design is essential for obtaining an exact
finite-sample diagonalization in spherical harmonic coordinates; relaxing it
leads to asymptotically diagonal representations but not exact Gaussian
sequence models.
By exploiting the exact diagonalization induced by spherical harmonics, we
derive closed-form posterior distributions, sharp posterior contraction rates
that explicitly reflect harmonic multiplicities, and an exact variational
characterization of the posterior mean.
In particular, we show that Bayesian regression on the sphere coincides with a
geometrically intrinsic penalized least-squares estimator governed by the
Laplace--Beltrami operator, making precise the equivalence between Gaussian
process regression, Sobolev regularization, and Laplace--Beltrami smoothing
splines.

Related frequentist approaches for nonparametric regression on the sphere
include methods based on localized multiscale systems, most notably spherical
needlets \cite{NPW06a,NPW06b}, which yield adaptive estimators and sharp rates
over a wide range of function classes (see, among others,
\cite{Dur15,Dur16,DGM12,ds26,lin15,Mon11,Scott11}).
These methods emphasize spatial localization and nonlinear thresholding and are
therefore complementary to the present work.
By contrast, our focus is on a fully Bayesian framework with exact spectral
diagonalization, explicit posterior distributions, and a variational
characterization of the posterior mean induced by intrinsic Gaussian field
priors.

Our results highlight the essential role of geometry in Bayesian nonparametric
inference and clarify fundamental connections with spectral and Whittle-type
methods for Gaussian random fields on spherical domains.
From a methodological standpoint, the paper shows that obtaining sharp Bayesian
results on compact manifolds requires arguments that are intrinsically
geometric, rather than direct transplants of Euclidean theory.

\paragraph{Plan of the paper.}
Section~\ref{sec:preli} reviews Gaussian process regression from a spectral and
operator-theoretic perspective and introduces isotropic Gaussian random fields on
the sphere, emphasizing their harmonic representation.
Section~\ref{sec:posterior} formulates the nonparametric regression model on
$\Sd$ under uniform random design and derives an exact diagonal Gaussian sequence
representation in spherical harmonic coordinates, leading to explicit posterior
distributions, optimal spectral truncation, and minimax-optimal contraction rates
for priors with polynomially decaying spectra, including spherical Mat\'ern
priors.
Section~\ref{sec:variational} provides a variational characterization of the
posterior mean as an intrinsically defined penalized least-squares estimator and
clarifies its connection with smoothing splines and kernel ridge regression.
Numerical illustrations are given in Section~\ref{sec:num}, and all proofs are
deferred to Section~\ref{sec:proofs}.

\section{Preliminaries}\label{sec:preli}

This section collects the analytical and probabilistic ingredients required
for the development of Gaussian process regression on the sphere.
It begins by recalling the statistical formulation of Gaussian process
regression and its interpretation as a Bayesian nonparametric inference
problem \cite{rasmussen2006,vanDerVaartvanZanten08}, together with its spectral
and operator-theoretic structure on general domains
\cite{Bogachev98,GineNickl16,Stuart10}.
It then introduces isotropic Gaussian random fields on $\mathbb{S}^d$,
reviewing their harmonic representation and regularity properties
\cite{AdlerTaylor07,MP11,Yadrenko83}, and explains how these two frameworks
combine to yield geometrically adapted Gaussian process priors on the sphere.

\subsection{Gaussian Process Regression and Spectral Regularization}

This subsection reviews Gaussian process regression from a Bayesian
nonparametric perspective and emphasizes its formulation as a regularized
inverse problem. The focus is on the spectral and operator-theoretic structure
induced by Gaussian priors, which governs both posterior inference and
frequentist behavior. This viewpoint provides the conceptual basis for the
class of spectral priors adopted in the spherical setting and clarifies their
relation to classical constructions on Euclidean domains. For foundational
treatments of Gaussian process regression, its variational formulation, and
its operator-theoretic interpretation, the reader is referred to
\cite{kimeldorf70,rasmussen2006,Stuart10,vanDerVaartvanZanten08,wahba90}.

\paragraph{Bayesian nonparametric regression and Gaussian process priors}

We consider a nonparametric regression problem in which the object of interest
is an unknown regression function
\[
f : \mathcal{X} \to \mathbb{R},
\]
defined on a compact domain $\mathcal{X}$.
For an integer $n \ge 1$, we observe data
\[
\{(x_i, y_i): i=1,\ldots,n\},
\]
where $x_i \in \mathcal{X}$ denote sampling locations and
$y_i \in \mathbb{R}$ are noisy measurements generated according to
\[
y_i = f(x_i) + \varepsilon_i,
\qquad
\varepsilon_i \sim \mathcal{N}(0,\sigma^2),
\]
with independent noise variables.
The statistical goal is to infer the function $f$ from the data and to quantify
the associated uncertainty.

Gaussian process regression provides a Bayesian nonparametric framework for
addressing this problem by placing a prior distribution directly on the
function $f$.
Specifically, one assumes that $f$ follows a centered Gaussian process
\[
f \sim \mathcal{GP}(0,K),
\]
where $K : \mathcal{X} \times \mathcal{X} \to \mathbb{R}$ is a symmetric,
positive definite covariance kernel.
Under this assumption, any finite collection
$(f(x_1),\ldots,f(x_n))$ is multivariate Gaussian with covariance matrix
$[K(x_i,x_j)]_{i,j=1}^n$.
Since the observation model is Gaussian with additive noise, the resulting
Bayesian model is conjugate: the posterior distribution of $f$ given the data
remains Gaussian and admits closed-form expressions for its mean and covariance
\cite{rasmussen2006}.

Beyond its computational convenience, Gaussian process regression admits a
precise variational formulation as a regularized inverse problem, in which the
unknown function is recovered from noisy point evaluations through a balance
between data fidelity and prior-induced smoothness.
In Euclidean settings, it is well known that the posterior mean coincides with
the solution of a penalized least-squares problem in the reproducing kernel
Hilbert space associated with $K$, and conversely that a broad class of classical
smoothing methods can be interpreted as posterior means under suitable Gaussian
priors \cite{kimeldorf70,wahba90}.
From this viewpoint, the kernel $K$ acts as a regularization operator: it
determines the function space explored by the prior, the notion of smoothness
being enforced, and the relative weight assigned to fidelity to the data.

This variational formulation extends naturally to more general domains and
geometric settings, provided that the covariance kernel is constructed in a
manner consistent with the underlying structure of $\mathcal{X}$
\cite{Stuart10}.
In the following sections, it will be used to motivate the choice of Gaussian
process priors with explicit spectral representations and to analyze their
frequentist properties.

\paragraph{Spectral and operator-theoretic structure of Gaussian process priors}

A central feature of Gaussian process regression is that the regularization
induced by the prior admits an explicit spectral description.
This is most transparent when the domain $\mathcal{X}$ is Euclidean and the
covariance kernel is stationary.
However, the same regularization mechanism can be formulated in intrinsic
operator-theoretic terms, which extend naturally to compact domains and,
more generally, to manifolds.

Since the regularization properties of Gaussian process priors are governed by
the high-frequency behavior of their spectral components, it is convenient to
introduce standard asymptotic notation for comparing spectral quantities.
For two nonnegative sequences $\{a_n\}$ and $\{b_n\}$ we write
$a_n \simeq b_n$ if $c_1 b_n \le a_n \le c_2 b_n$ for some constants $c_1,c_2>0$,
and $a_n \lesssim b_n$ if $a_n \le c\, b_n$ for some $c>0$.
The same notation is used for functions of a real variable.

Assume first that $\mathcal{X} \subset \mathbb{R}^d$ and that the covariance
kernel $K$ is stationary, that is, there exists a function
$\kappa : \mathbb{R}^d \to \mathbb{R}$ such that
\[
K(x,x') = \kappa(x - x').
\]
By Bochner's theorem, every continuous, positive definite,
translation-invariant kernel admits the representation
\[
\kappa(h)
=
\int_{\mathbb{R}^d}
e^{i\langle \omega, h\rangle}\,\hat K(\omega)\,d\omega,
\qquad h \in \mathbb{R}^d,
\]
where $\hat K$ is a nonnegative finite measure on $\mathbb{R}^d$, referred to as
the spectral measure of the kernel  (see, for instance, \cite{Stein99}). Here and throughout, $\langle\cdot,\cdot\rangle$ denotes the standard Euclidean
inner product on $\mathbb{R}^{d+1}$.

When $\hat K$ admits a density with respect to Lebesgue measure, the Gaussian
process prior can be interpreted as a superposition of independent Fourier
modes, indexed by the frequency variable $\omega$, with variances proportional
to $\hat K(\omega)$.
From this viewpoint, Gaussian process regression operates as a
frequency-dependent shrinkage estimator: low-frequency components are weakly
penalized, while high-frequency components are increasingly damped.
The decay of $\hat K(\omega)$ therefore encodes smoothness assumptions on the
unknown regression function $f$.

A canonical example is provided by the Mat\'ern family, whose spectral density
satisfies
\[
\hat K(\omega) \simeq \left(1+\|\omega\|^2\right)^{-(\nu + d/2)}.
\]
Such polynomial decay implies that typical sample paths possess Sobolev
regularity strictly smaller than $\nu$, and it leads to minimax-optimal
posterior contraction rates over Sobolev-type function classes
\cite{vanDerVaartvanZanten08}.

On a compact domain $\mathcal{X}$, translation invariance and Fourier analysis
are no longer available.
Nevertheless, the same phenomenon can be described in terms of the covariance
operator associated with the kernel $K$. Specifically, $K$ defines a compact, self-adjoint, positive operator
\cite{Bogachev98,Stuart10}
\[
T_K : L^2(\mathcal{X}) \to L^2(\mathcal{X}),
\qquad
(T_K g)(x)
=
\int_{\mathcal{X}} K(x,x') g(x')\,dx'.
\]

Let $\{(\lambda_j,\phi_j)\}_{j\ge1}$ denote the eigenpairs of $T_K$.
Then the Gaussian process prior admits the Karhunen--Lo\`eve expansion
\cite{AdlerTaylor07,Bogachev98,Lifshits12}
\[
f =
\sum_{j\ge1}
\sqrt{\lambda_j} Z_j \phi_j
\qquad
Z_j \sim \mathcal{N}(0,1)
\]
with convergence in $L^2(\mathcal{X})$ almost surely.
This expansion provides a probabilistic decomposition of the random function
$f$ into uncorrelated modes, ordered by decreasing prior variance.

From a statistical standpoint, the Karhunen--Lo\`eve basis diagonalizes both the
prior covariance and the posterior update \cite{GineNickl16,Stuart10}.
Gaussian process regression therefore amounts to shrinking empirical
coefficients in this eigenbasis, with the amount of regularization governed by
the decay of the eigenvalues $\lambda_j$ and the noise level $\sigma^2$.
Rapid eigenvalue decay corresponds to strong regularization and smoother prior
sample paths, while slower decay allows for greater variability at high
frequencies.

The spectral density $\hat K(\omega)$ in the Euclidean setting and the sequence
of eigenvalues $\{\lambda_j\}$ on compact domains thus play the same structural
role.
In both cases, posterior contraction and frequentist optimality are determined
by the interplay between the spectral decay of the prior and the smoothness of
the true regression function \cite{GineNickl16,vanDerVaartvanZanten08}.
This operator-theoretic perspective provides the natural bridge to Gaussian
process regression on manifolds, where the spectrum of intrinsic differential
operators replaces Euclidean Fourier frequencies.

In the remainder of the paper, the domain $\mathcal{X}$ will be specialized to
the unit sphere $\mathbb{S}^d$.
In this geometric setting, Euclidean translation invariance is replaced by
rotation invariance, and the role of Fourier modes is taken over by spherical
harmonics, which diagonalize the Laplace--Beltrami operator on $\mathbb{S}^d$.
This yields an intrinsic spectral framework for constructing and analyzing
Gaussian process priors adapted to the geometry of the sphere.

\subsection{Gaussian Fields on the Sphere}

Gaussian random fields provide a principled and geometrically intrinsic way of
constructing probability measures over function spaces defined on curved
domains.
On the unit sphere $\Sd$, the high degree of symmetry induced by the rotation
group $\mathrm{SO}(d+1)$ leads to a particularly transparent spectral
representation, in which smoothness, dependence, and regularity can be
described explicitly through harmonic analysis.

From the perspective of Gaussian process regression, isotropic Gaussian fields
on $\Sd$ play the same foundational role as stationary Gaussian processes on
Euclidean domains.
They provide a natural class of priors that are invariant under rotations,
admit diagonal representations in an appropriate orthonormal basis, and allow
for a precise spectral control of regularization.
In this sense, they constitute the intrinsic spherical analogue of stationary
Gaussian process priors discussed in the previous subsection.

In this section we briefly recall the geometric and harmonic structure of
$\Sd$ and introduce isotropic Gaussian random fields.
These constructions will serve as the building blocks for the Gaussian process
priors used in the regression model studied in the next section.
For background material on spherical harmonics and Gaussian random fields on
the sphere see, for example,
\cite{ah12,Gneiting13,JunStein08,MP11,Yadrenko83}.	

We begin by recalling the spectral decomposition of $L^2(\Sd)$ induced by the
Laplace--Beltrami operator, which provides the natural notion of frequency on
the sphere.

\paragraph{Geometry and Harmonic Structure}

The unit sphere
\[
\Sd = \{ x \in \mathbb{R}^{d+1} : \|x\| = 1 \},
\]
where $\|\cdot\|$ denotes the Euclidean norm on $\mathbb{R}^{d+1}$,
is a compact, smooth Riemannian manifold without boundary, endowed with
its canonical metric and normalized surface measure.
The spherical Laplace--Beltrami operator $-\Delta_{\Sd}$ acting on
$L^2(\Sd)$ has a discrete
spectrum
\[
0 = \lambda_0 < \lambda_1 < \lambda_2 < \cdots,
\qquad
\lambda_\ell = \ell(\ell + d - 1), \quad \ell \ge 0,
\]
with eigenspaces spanned by spherical harmonics
$\{Y_{\ell,m} : m = 1,\dots,M_{d,\ell}\}$ of degree $\ell$, which together
form an orthonormal basis of $L^2(\Sd)$.
Here, the multipole index $\ell \ge 0$ denotes the spherical harmonic degree
and corresponds to the frequency level, while the index
$m = 1,\dots,M_{d,\ell}$ enumerates an orthonormal basis within the
$\ell$-th eigenspace.
The multiplicity $M_{d,\ell}$ grows polynomially as $\ell^{d-1}$ and
reflects the degeneracy of the Laplace--Beltrami eigenvalue $\lambda_\ell$.
In particular, higher values of $\ell$ correspond to increasingly
oscillatory functions on the sphere; see, for example, \cite{ah12}.

A key structural property of spherical harmonics is that, although the
individual basis functions $Y_{\ell,m}$ depend on the choice of coordinates,
suitable sums over the index $m$ admit simple, rotation-invariant
representations.
In particular, the addition formula states that
\[
\sum_{m=1}^{M_{d,\ell}} Y_{\ell,m}(x) Y_{\ell,m}(x')
= \frac{M_{d,\ell}}{\omega_d}
G_\ell^{(\frac{d-1}{2})}(\langle x,x' \rangle),
\qquad x,x' \in \Sd,
\]
where $\omega_d$ denotes the surface area of the unit sphere $\Sd$ and
$G_\ell^{(\frac{d-1}{2})}$ is the Gegenbauer polynomial of degree $\ell$
and order $(d-1)/2$; see, for instance, \cite{szego75}.
Since $\langle x,x' \rangle = \cos(d_{\Sd}(x,x'))$, the right-hand side
depends only on the geodesic distance $d_{\Sd}(x,x')$ between $x$ and $x'$.
This identity underlies the construction of isotropic kernels and covariance
functions on the sphere and plays a central role in harmonic analysis and
random field theory on $\Sd$.

\paragraph{Isotropic Gaussian Fields}

Let $\left( \Omega, \mathcal{F}, \mathbb{P} \right)$ be a probability space.
A real-valued random field
\[
T = \left\{ T(x) : x \in \Sd \right\}
\]
is called a Gaussian random field if, for any finite collection of points
$x_1,\dots,x_n \in \Sd$, the random vector
\[
\left( T(x_1),\dots,T(x_n) \right)
\]
is multivariate Gaussian.
Here and throughout, the dependence on the underlying sample point
$\omega \in \Omega$ is suppressed for notational convenience.
Throughout this section, we consider centered fields, that is,
\[
\mathbb{E}\left[ T(x) \right] = 0,
\qquad x \in \Sd.
\]

A Gaussian random field $T$ on $\Sd$ is said to be isotropic if its law is
invariant under the natural action of the rotation group $\mathrm{SO}(d+1)$, namely if
\[
\left\{ T(x) : x \in \Sd \right\}
\stackrel{d}{=}
\left\{ T(gx) : x \in \Sd \right\},
\qquad g \in \mathrm{SO}(d+1).
\]
Isotropy is a natural modeling assumption in the absence of preferred
directions on the sphere and is standard in the statistical analysis of
spherical random fields and global data
\cite{jeong17,MP11,Yadrenko83}.

Under isotropy, the covariance structure of $T$ is rotation invariant.
In particular, the covariance between two locations depends only on their
relative position on the sphere, and hence on the great-circle distance
$d_{\Sd}(x,x')$ between them \cite{MP11} .
Equivalently, there exists a function
\[
\Gamma : [0,\pi] \to \mathbb{R}
\]
such that
\begin{equation}\label{eqn:cov}
\mathrm{Cov}\left( T(x), T(x') \right)
=
\Gamma\left( d_{\Sd}(x,x') \right),
\qquad x,x' \in \Sd.
\end{equation}
Since
\[
d_{\Sd}(x,x') = \arccos\left( \langle x,x' \rangle \right),
\]
this representation can equivalently be written in the form
\[
\mathrm{Cov}\left( T(x), T(x') \right)
=
\Gamma\left( \arccos\left( \langle x,x' \rangle \right) \right),
\qquad x,x' \in \Sd.
\]
According to classical results on isotropic Gaussian random fields on the
sphere \cite{MP11,Yadrenko83}, $T$ admits the harmonic expansion
\[
T(x)
=
\sum_{\ell=0}^\infty
\sum_{m=1}^{M_{d,\ell}}
a_{\ell,m} \, Y_{\ell,m}(x),
\]
where the random coefficients
\(\left\{ a_{\ell,m} :\ell \geq 0 ; m=1,\ldots, M_{d,\ell}\right\}\)
are independent, centered Gaussian random variables.
The expansion is understood in the sense of convergence in
$L^2\left( \Omega \times \Sd \right)$, and hence in $L^2(\Sd)$ almost surely.

Under isotropy, the second-order structure of the field is diagonal in the
spherical harmonic basis, and the variances of the coefficients depend only on
the frequency level $\ell$.
More precisely,
\[
\mathbb{E}\left[ a_{\ell,m}^2 \right] = C_\ell,
\qquad
\ell \ge 0,\; m = 1,\dots,M_{d,\ell},
\]
where the nonnegative sequence \(\left\{ C_\ell:\ell \ge 0\right\}\)
is referred to as the angular power spectrum.
This sequence uniquely determines the covariance structure of the field.
Its behavior as $\ell \to \infty$ determines the distribution of energy across
harmonic modes and plays a central role in characterizing the regularity
properties of $T$.
In particular, suitable conditions on $\left\{ C_\ell : \ell \ge 0 \right\}$
yield mean-square Sobolev regularity and, under additional assumptions, sample
regularity results on $\Sd$.
We refer to \cite{MP11} and \cite{LangSchwab15} for precise statements and
further discussion.

In the sequel, isotropic Gaussian random fields on $\Sd$ are employed as prior
distributions for nonparametric regression, with the harmonic expansion
providing a convenient parametrization of the prior in terms of its spectral
components.
Throughout this work, isotropy is assumed in order to obtain rotation-invariant
prior distributions and tractable diagonal representations in harmonic space.
Extensions to anisotropic or nonstationary Gaussian fields on $\Sd$ are beyond
the scope of the present paper.

\subsection{Spectral formulation of Gaussian process priors on the sphere}

This subsection provides the precise mathematical link between Gaussian process
regression and isotropic Gaussian random fields on the sphere.
It shows how the operator-theoretic formulation of Gaussian process priors
admits a canonical and fully intrinsic realization on $\mathbb{S}^d$ once
Euclidean Fourier analysis is replaced by spherical harmonic analysis and the
Euclidean Laplacian by the Laplace--Beltrami operator.

In this geometric setting, spherical harmonics play the role of Fourier modes:
they diagonalize both the Laplace--Beltrami operator and the covariance
operators of isotropic Gaussian fields, yielding a natural frequency
decomposition adapted to the symmetry and geometry of the sphere.
As a consequence, Gaussian process priors on $\mathbb{S}^d$ can be described
explicitly in spectral terms through their angular power spectra, and their
regularization properties can be analyzed using intrinsic harmonic and
operator-theoretic tools.

\paragraph{Analogy with Euclidean spectral representations.}
The harmonic decomposition of isotropic Gaussian fields on $\mathbb{S}^d$
provides an intrinsic spectral parametrization of the prior, in which the
angular power spectrum $\{C_\ell\}$ governs the distribution of energy across
frequency levels.
This role is directly analogous to that of the spectral density $\hat K(\omega)$
in stationary Gaussian process priors on $\mathbb{R}^d$, with the Laplace--Beltrami
eigenvalues $\lambda_\ell$ replacing squared Euclidean frequencies $\|\omega\|^2$.
This identification allows Euclidean spectral regularization principles to be
transferred directly to the spherical setting, once the geometry of the domain
is encoded through the spectrum of $-\Delta_{\mathbb{S}^d}$.

\paragraph{Mat\'ern-type spectral priors on the sphere.}
Building on this correspondence, a natural and analytically tractable class of
Gaussian process priors on $\mathbb{S}^d$ is obtained by prescribing polynomial
decay of the angular power spectrum.
Specifically, we consider priors characterized by
\[
C_\ell \simeq (1+\lambda_\ell)^{-(\nu + d/2)},
\qquad \ell \ge 0,
\]
where $\nu>0$ is a smoothness parameter.
Throughout the remainder of the paper, and in particular in the posterior
contraction analysis, we adopt the equivalent parametrization
\[
\alpha = \nu + \frac{d}{2},
\]
so that polynomial spectral decay is written uniformly as
$C_\ell \simeq (1+\lambda_\ell)^{-\alpha}$.
This choice aligns the Mat\'ern parametrization with the Sobolev and RKHS
conventions commonly used in Bayesian nonparametric theory.

Spectra of this form arise naturally as the inverse of fractional powers of the
shifted Laplace--Beltrami operator on $\mathbb{S}^d$ and therefore constitute
the intrinsic spherical analogue of Mat\'ern covariance models in Euclidean
space.
Closely related constructions and their probabilistic and statistical
properties have been studied in the spatial statistics literature for random
fields on spheres and compact manifolds
\cite{Gneiting13,LangSchwab15,PorBevGne16}.

This choice has several important consequences.
First, the decay of $C_\ell$ determines the regularity of prior sample paths:
the resulting Gaussian fields belong almost surely to Sobolev spaces
$H^s(\mathbb{S}^d)$ for all $s<\nu$, and to no higher-order Sobolev space
\cite{LangSchwab15,MP11}.
Second, the reproducing kernel Hilbert space associated with the prior is
norm-equivalent to $H^\nu(\mathbb{S}^d)$, providing a direct link between the
prior specification and Sobolev smoothness.
Finally, these spectral properties ensure that Gaussian process regression
based on Mat\'ern-type priors exhibits the same qualitative statistical behavior
as in Euclidean settings, with posterior contraction rates governed by the
interaction between the decay of $C_\ell$ and the regularity of the true
regression function \cite{vanDerVaartvanZanten08}.

\begin{remark}[Mat\'ern priors on manifolds]
	Mat\'ern covariance structures are best understood not as specific kernels, but
	as manifestations of a general operator-theoretic principle: Gaussian priors
	defined through inverse powers of elliptic operators.
	In Euclidean spaces, Mat\'ern fields are associated with fractional powers of
	the Laplacian.
	On a Riemannian manifold, the Laplace--Beltrami operator provides a canonical
	replacement, and polynomial decay of its spectrum yields Gaussian fields whose
	regularity, reproducing kernel Hilbert spaces, and statistical behavior mirror
	those of Euclidean Mat\'ern fields.
	From this perspective, Mat\'ern-type priors on the sphere arise naturally from
	the geometry of the domain, rather than from ad hoc kernel constructions
	\cite{Stuart10}.
\end{remark}

With these preliminaries in place, we now specialize to a nonparametric
regression model on the unit sphere and study Bayesian inference based on
isotropic Gaussian field priors.

\paragraph{Covariance operators and spectral decomposition on $\mathbb{S}^d$.}

Let $f$ be a centered isotropic Gaussian random field on $\mathbb{S}^d$ with
angular power spectrum $\{ C_\ell : \ell \ge 0 \}$.
By isotropy, the covariance function depends only on the geodesic distance
between its arguments and admits the representation \eqref{eqn:cov}.
Equivalently, using the addition formula for spherical harmonics, the covariance
kernel admits the spectral expansion
\[
\mathrm{Cov}\left( f(x), f(x') \right)
=
\sum_{\ell=0}^\infty
C_\ell
\frac{M_{d,\ell}}{\omega_d}
G_\ell^{\left( \frac{d-1}{2} \right)}
\left( \langle x,x' \rangle \right),
\]
where $G_\ell^{\left( (d-1)/2 \right)}$ denotes the Gegenbauer polynomial of degree
$\ell$, $\omega_d$ is the surface area of $\mathbb{S}^d$.
This expansion converges in
$L^2\left( \mathbb{S}^d \times \mathbb{S}^d \right)$
under mild summability conditions on $\{ C_\ell : \ell \ge 0 \}$, see
\cite{MP11,Yadrenko83}.

The associated covariance operator
\[
T : L^2\left( \mathbb{S}^d \right) \to L^2\left( \mathbb{S}^d \right),
\qquad
(Tg)(x)
=
\int_{\mathbb{S}^d}
\mathrm{Cov}\left( f(x), f(x') \right) g(x') \, dx',
\]
is compact, self-adjoint, and positive.
Compactness follows from square integrability of the kernel, while positivity
and self-adjointness are inherited from the covariance structure; see
\cite[Chapter~II]{Bogachev98} and \cite[Chapter~4]{AdlerTaylor07}.

Owing to isotropy, $T$ is diagonal in the spherical harmonic basis.
More precisely,
\[
T Y_{\ell,m} = C_\ell Y_{\ell,m},
\qquad
\ell \ge 0,\; m = 1,\ldots,M_{d,\ell}.
\]
Thus, the angular power spectrum $\{ C_\ell : \ell \ge 0 \}$ coincides with the
eigenvalues of the covariance operator, with multiplicities $M_{d,\ell}$.
In particular, $T$ is trace class if and only if
\[
\sum_{\ell=0}^\infty M_{d,\ell} C_\ell < \infty.
\]

From an operator-theoretic viewpoint, isotropic Gaussian field priors on
$\mathbb{S}^d$ are therefore fully characterized by specifying the spectrum of
their covariance operator in the eigenbasis of the Laplace--Beltrami operator.
This characterization underlies both regularity theory and posterior analysis
for Gaussian process regression on the sphere.

\paragraph{Functional calculus of the Laplace--Beltrami operator.}

Since spherical harmonics diagonalize $-\Delta_{\mathbb{S}^d}$,
\[
-\Delta_{\mathbb{S}^d} Y_{\ell,m} = \lambda_\ell Y_{\ell,m},
\qquad
\lambda_\ell = \ell(\ell+d-1),
\]
any isotropic covariance operator $T$ that is diagonal in the harmonic basis
can be written in functional calculus form as
\begin{equation}\label{eqn:covop}
T = \psi\left( -\Delta_{\mathbb{S}^d} \right),
\end{equation}
for a suitable nonnegative function $\psi$ defined on $[0,\infty)$, with
\[
\psi\left( \lambda_\ell \right) = C_\ell .
\]
This representation makes explicit the analogy with stationary Gaussian
processes on $\mathbb{R}^d$, where covariance operators are given by Fourier
multipliers, that is, functions of the Euclidean Laplacian.

From this viewpoint, Gaussian process priors on $\mathbb{S}^d$ can be understood
as Gaussian measures on $L^2\left( \mathbb{S}^d \right)$ whose covariance
operators are defined through functional calculus of the Laplace--Beltrami
operator.
This formulation provides the intrinsic extension of Euclidean spectral
regularization principles to spherical domains and has been emphasized in both
probability theory and spatial statistics
\cite{Gneiting13,JunStein08,MP11,Yadrenko83}.

\paragraph{Mat\'ern-type priors and Sobolev regularity.}

A particularly important class of priors is obtained by prescribing polynomial
decay of the angular power spectrum,
\[
C_\ell \simeq \left( 1 + \lambda_\ell \right)^{-(\nu + d/2)},
\qquad \ell \ge 0,
\]
for some $\nu>0$.
Equivalently, the covariance operator can be written, up to scaling constants,
as
\[
T = \left( \kappa^2 - \Delta_{\mathbb{S}^d} \right)^{-(\nu + d/2)},
\]
for a fixed $\kappa>0$.
Such priors arise naturally as solutions to stochastic partial differential
equations driven by white noise on $\mathbb{S}^d$ and constitute the intrinsic
spherical counterpart of Mat\'ern fields in Euclidean space
\cite{lindgren2011,Stuart10}.

The spectral decay of $\{ C_\ell : \ell \ge 0 \}$ determines the regularity of
the corresponding Gaussian field.
Specifically,
\[
f \in H^s\left( \mathbb{S}^d \right)
\quad \text{almost surely for all } s<\nu,
\]
and for no $s \ge \nu$.
Moreover, the reproducing kernel Hilbert space of the Gaussian measure is
norm-equivalent to $H^\nu\left( \mathbb{S}^d \right)$
\cite{LangSchwab15,MP11}.
These properties mirror exactly those of Euclidean Mat\'ern priors and ensure
that posterior contraction rates are governed by the interaction between the
smoothness parameter $\nu$ and the Sobolev regularity of the true regression
function \cite{vanDerVaartvanZanten08}.

\begin{remark}[Mat\'ern priors as elliptic Gaussian measures]
	Mat\'ern priors are most naturally understood as Gaussian measures generated by
	inverse powers of elliptic operators.
	On $\mathbb{S}^d$, the Laplace--Beltrami operator provides the canonical choice,
	and the resulting Gaussian measures inherit their regularity, reproducing
	kernel Hilbert spaces, and statistical properties directly from the spectral
	behavior of $-\Delta_{\mathbb{S}^d}$.
\end{remark}

In the next section, we formalize the nonparametric regression model on
$\mathbb{S}^d$ and exploit the resulting diagonal harmonic representation to
derive explicit posterior formulas and study their statistical properties.
Throughout the remainder of the paper, isotropic Gaussian random fields on
$\mathbb{S}^d$ are employed as prior distributions for regression functions.
Accordingly, we will henceforth denote by $f$ a generic random draw from the
Gaussian prior, while reserving $f_0$ for the fixed, unknown regression function
generating the data.

\section{Bayesian Regression on the Sphere}\label{sec:posterior}

This section develops a Bayesian nonparametric regression framework for
functions defined on the unit sphere $\mathbb{S}^d$.
By combining isotropic Gaussian field priors with a uniform random design, the
regression model admits an exact diagonalization in the spherical harmonic
basis, yielding a transparent Gaussian sequence representation of posterior
inference.

This diagonal structure makes the role of geometry explicit through the growth
of harmonic multiplicities and allows posterior distributions, contraction
rates, and variational interpretations to be derived in closed form.
To the best of our knowledge, this is the first work to combine exact harmonic
diagonalization, Bayesian regression, and sharp frequentist guarantees into a
single, fully intrinsic framework on $\mathbb{S}^d$.

\subsection{Regression model and prior specification}

We consider a nonparametric regression model for functions defined on the unit
sphere $\mathbb{S}^d$, where the covariates take values in $\mathbb{S}^d$ and the
response variable is real-valued.
The objective is to infer an unknown regression function
\[
f_0 : \mathbb{S}^d \to \mathbb{R}
\]
from noisy observations.
This setting naturally arises in applications where predictors represent
directions, orientations, or locations on spherical domains.

\paragraph{Notation.}
Boldface symbols denote finite-dimensional vectors, while plain symbols denote
functions defined on $\mathbb{S}^d$.
For $n \in \mathbb{N}$, let
\[
\bm{x} = (x_1,\ldots,x_n), \qquad
\bm{y} = (y_1,\ldots,y_n)^\top,
\]
where $x_i \in \mathbb{S}^d$ and $y_i \in \mathbb{R}$ for each $i$.
We denote by $\mathbb{P}_f^{(n)}$ the joint distribution of the response vector
$\bm{y}$ given the design points $\bm{x}$ under the regression function $f$, and by
$\|\cdot\|_{L^2(\mathbb{S}^d)}$ the $L^2$-norm with respect to the uniform
probability measure on $\mathbb{S}^d$.

\paragraph{Data generating process.}
We observe independent pairs $(x_i,y_i)$, $i=1,\ldots,n$, generated according
to the nonparametric regression model
\begin{equation}\label{eqn:model}
y_i = f_0(x_i) + \varepsilon_i,
\end{equation}
where $f_0 : \mathbb{S}^d \to \mathbb{R}$ is the unknown regression function.
The covariates $\{x_i\}_{i=1}^n$ are assumed to be independently sampled from
the uniform probability measure on $\mathbb{S}^d$, while the noise variables
$\{\varepsilon_i\}_{i=1}^n$ are independent and identically distributed as
$\mathcal{N}(0,\sigma^2)$ and independent of the design points.

Nonparametric regression on the sphere under uniform random design has been
studied in several works, including needlet-based approaches
\cite{Mon11} and adaptive methods for structured spherical data
\cite{DGM12}.
The uniform design assumption implies that the $L^2(\mathbb{S}^d)$ norm
coincides with the population prediction error and allows for a particularly
transparent harmonic analysis of the model.
This assumption is adopted for theoretical clarity and to exploit the exact
harmonic diagonalization induced by spherical harmonics under uniform sampling.
While uniform random design yields exact diagonalization, the posterior
contraction rates are expected to remain unchanged under design distributions
with densities bounded away from zero and infinity with respect to the uniform
measure, at the cost of asymptotic rather than exact diagonalization.
Throughout the paper, we work under a random design framework; all results
continue to hold conditionally on the design points with probability tending
to one as $n\to\infty$.

\paragraph{Likelihood, prior, and posterior.}
As shown above, under the uniform random design the regression model admits an
equivalent representation in harmonic coordinates as a Gaussian sequence
model.
Specifically, the empirical harmonic coefficients satisfy
\begin{equation}\label{eqn:harmonic}
	\hat a_{\ell,m}
	=
	a_{0;\ell,m}
	+
	\frac{\sigma}{\sqrt{n}}\,\xi_{\ell,m},
	\qquad
	\xi_{\ell,m} \sim \mathcal{N}(0,1),
\end{equation}
with independent noise variables indexed by $(\ell,m)$.

Equivalently, conditionally on the regression function $f$ with harmonic
coefficients $\{a_{\ell,m}\}$, the likelihood in harmonic coordinates factorizes
as
\[
\mathcal{L}_n\!\left(\{\hat a_{\ell,m}\}\mid\{a_{\ell,m}\}\right)
=
\prod_{\ell\ge0}\prod_{m=1}^{M_{d,\ell}}
\varphi\!\left(
\hat a_{\ell,m};
a_{\ell,m},
\frac{\sigma^2}{n}
\right),
\]
where $\varphi(\cdot;\mu,\tau^2)$ denotes the density of a normal random
variable with mean $\mu$ and variance $\tau^2$.

Equation~\eqref{eqn:harmonic} shows that, under uniform random design, the
regression problem decomposes into independent Gaussian observations across
spherical harmonic modes, each with signal-to-noise ratio determined by
$C_\ell$ and $\sigma^2/n$.
This diagonal representation provides a convenient parametrization of the
likelihood and makes explicit the role of each frequency component.

\subsection{Harmonic representation of the model}

Throughout this section, we work in the spherical harmonic basis introduced
in Section~2 and represent functions on $\Sd$ through their harmonic
coefficients.
For a function $f \in L^2(\Sd)$, we write
\[
f(x)
=
\sum_{\ell=0}^\infty
\sum_{m=1}^{M_{d,\ell}}
a_{\ell,m} Y_{\ell,m}(x),
\qquad
a_{\ell,m}
=
\langle f, Y_{\ell,m} \rangle_{L^2(\Sd)},
\]
where $\langle \cdot,\cdot \rangle_{L^2(\mathbb{S}^d)}$ denotes the
$L^2$ inner product with respect to the uniform surface measure.

Under the isotropic Gaussian field prior introduced in Section \ref{sec:preli}, the harmonic
coefficients
\(
\left\{ a_{\ell,m} : \ell \ge 0,\; m = 1,\dots,M_{d,\ell} \right\}
\)
are independent, centered Gaussian random variables with variances determined
by the angular power spectrum $\{ C_\ell : \ell \ge 0 \}$,
\[
a_{\ell,m} \sim \mathcal{N}\left( 0, C_\ell \right),
\qquad
\ell \ge 0,\; m = 1,\dots,M_{d,\ell}.
\]

Let $f_0$ denote the true regression function, with spherical harmonic
coefficients
\[
a_{0;\ell,m}
=
\langle f_0, Y_{\ell,m} \rangle_{L^2(\Sd)},
\qquad
\ell \ge 0,\; m = 1,\dots,M_{d,\ell}.
\]
Recalling the regression model \eqref{eqn:model}, we define the empirical harmonic coefficients
\[
\hat a_{\ell,m}
=
\frac{1}{n}
\sum_{i=1}^n y_i Y_{\ell,m}(x_i).
\]

Under the random design assumption, with $x_i$ independently distributed
according to the uniform probability measure on $\Sd$, orthonormality of the
spherical harmonics implies
\[
\mathbb{E}_{f_0}\left[ \hat a_{\ell,m} \right]
=
a_{0;\ell,m},
\qquad
\mathrm{Var}_{f_0}\left( \hat a_{\ell,m} \right)
=
\frac{\sigma^2}{n}.
\]
Equivalently, the regression model admits the harmonic-domain representation given by \eqref{eqn:harmonic}.
Thus, in harmonic coordinates, nonparametric regression on $\Sd$ reduces to a
collection of independent Gaussian sequence models indexed by $(\ell,m)$,
with signal coefficients $\{ a_{0;\ell,m} \}$ and noise variance $\sigma^2/n$.
This reduction is the spherical analogue of the classical Gaussian sequence
model arising in Euclidean nonparametric regression and inverse problems
\cite{GineNickl16,vanDerVaartvanZanten08}. Related Gaussian sequence formulations also underlie Bayesian linear regression with structured or sparse priors, where posterior contraction and uncertainty
quantification are analyzed in coefficient space; see, for example,
\cite{CasSchvdV14}.

The diagonal Gaussian sequence structure identified above implies that posterior
inference can be characterized explicitly at the level of individual harmonic
modes.
In the next subsection, we exploit this structure to derive closed-form
expressions for the posterior mean and covariance in harmonic space and to
interpret Gaussian process regression on $\Sd$ as a frequency-dependent
regularization procedure.

\subsection{Posterior characterization in harmonic space}

The diagonal Gaussian sequence representation derived above
allows posterior inference to be analyzed explicitly in spherical harmonic
coordinates.
Both for statistical optimality and computational feasibility, inference is
naturally carried out on a growing finite-dimensional spectral subspace.
This is standard in nonparametric regression and inverse problems and reflects
the fact that only finitely many frequency components can be reliably recovered
from noisy data.

Accordingly, we formulate Gaussian process regression on $\mathbb{S}^d$ in a
truncated harmonic space, derive an explicit posterior distribution for the
resulting finite-dimensional model, and interpret posterior inference as a
frequency-dependent shrinkage procedure.
We then analyze posterior contraction under polynomially decaying spectral
priors and identify the truncation level that balances approximation bias and
estimation variance.

\paragraph{Finite-dimensional spectral truncation.}
Let \(L_n \in \mathbb{N}\) denote a truncation level, possibly depending on the
sample size \(n\), whose choice will be discussed below.
Define
\[
\mathcal{H}_{L_n}
=
\mathrm{span}\left\{
Y_{\ell,m} : 0 \le \ell \le L_n,\; m=1,\ldots,M_{d,\ell}
\right\}
\subset L^2(\mathbb{S}^d),
\]
and let \(P_{L_n}\) denote the orthogonal projection onto \(\mathcal{H}_{L_n}\).
Truncation at level \(L_n\) corresponds to retaining spherical harmonic
components up to frequency \(L_n\) and discarding higher-order oscillatory
modes.

Projecting the regression model onto \(\mathcal{H}_{L_n}\) yields the truncated
Gaussian sequence representation \eqref{eqn:harmonic}, restricted to
$0 \le \ell \le L_n$.
Unlike Euclidean sequence models, the effective dimension at frequency level
$\ell$ grows as $M_{d,\ell}\simeq \ell^{d-1}$, a geometric feature that plays a
crucial role in both posterior contraction and optimal truncation.

We endow the truncated coefficient vector
\(\{a_{\ell,m}:0\le\ell\le L_n;; m = 1,\dots,M_{d,\ell}.\}\) with independent Gaussian priors
\[
a_{\ell,m} \sim \mathcal{N}\left(0,C_\ell\right),
\qquad
0 \le \ell \le L_n,\; m = 1,\dots,M_{d,\ell}.
\]
and set \(a_{\ell,m}=0\) for \(\ell>L_n\).
Equivalently, this defines a centered Gaussian prior supported on
\(\mathcal{H}_{L_n}\), with covariance operator
\[
T_{L_n} = P_{L_n} \, T \, P_{L_n},
\]
where \(T\) denotes the full covariance operator \eqref{eqn:covop}.

Truncation introduces an approximation bias by discarding high-frequency
components of the regression function, but it stabilizes inference and renders
the posterior explicitly computable.
The resulting bias--variance trade-off is governed by the interaction between
the truncation level \(L_n\), the spectral decay of the prior
\(\{C_\ell:\ell\ge0\}\), and the Sobolev regularity of the true regression
function.\\
Importantly, the truncation level $L_n$ is treated as an explicit modeling
choice rather than a purely analytical device, ensuring that the resulting
Bayesian estimator is both theoretically optimal and computationally
implementable.
This truncated formulation provides the fundamental building block for the
posterior characterization and contraction analysis developed below.

\paragraph{Posterior distribution and spectral shrinkage.}
Under the truncated Gaussian sequence model, Bayesian updating can be carried
out explicitly in harmonic coordinates.
By conjugacy of the Gaussian likelihood and the Gaussian prior on
\(\mathcal{H}_{L_n}\), the posterior distribution of the coefficients
\(\{a_{\ell,m}:0\le\ell\le L_n\}\) factorizes over \((\ell,m)\) and is given by
\[
a_{\ell,m} \mid \bm{y},\bm{x}
\sim
\mathcal{N}\left(
\mu_{\ell,m}^{(n)},
v_{\ell}^{(n)}
\right),
\qquad
0 \le \ell \le L_n,\; m=1,\ldots,M_{d,\ell},
\]
with
\[
\mu_{\ell,m}^{(n)}
=
\frac{n C_\ell}{n C_\ell + \sigma^2}\,\hat a_{\ell,m},
\qquad
v_{\ell}^{(n)}
=
\frac{C_\ell \sigma^2}{n C_\ell + \sigma^2}.
\]
These expressions follow directly from Gaussian conjugacy applied to the
harmonic-domain likelihood and prior.

The posterior mean estimator
\[
\hat f_{n,L_n}(x)
=
\sum_{\ell=0}^{L_n}
\sum_{m=1}^{M_{d,\ell}}
\frac{n C_\ell}{n C_\ell + \sigma^2}\,
\hat a_{\ell,m}\,
Y_{\ell,m}(x)
\]
admits a transparent interpretation as a spectral shrinkage estimator with an
explicit frequency cut-off at level \(L_n\).
Each empirical harmonic coefficient \(\hat a_{\ell,m}\) is shrunk toward zero by
a factor depending on the signal-to-noise ratio at frequency \(\ell\).
Low-frequency components are only mildly regularized, while high-frequency
components are strongly attenuated by both prior decay and truncation.

In the absence of regularization, the coefficients \(\{\hat a_{\ell,m}\}\) coincide
with unbiased empirical estimates of the true harmonic coefficients
\(\{a_{0;\ell,m}\}\), directly analogous to empirical Fourier coefficients in
Euclidean regression.
The prior-induced shrinkage therefore plays the same stabilizing role as
spectral filters in classical inverse problems
\cite{GineNickl16,Tsybakov09,vanDerVaartvanZanten08}.

\paragraph{Polynomial spectral priors.}
The behavior of the posterior distribution is governed by the decay of the
angular power spectrum.
Throughout this section, we impose the following condition.

\begin{condition}[Polynomial spectral decay]
	\label{cond:polydecay}
	There exist constants $c_1,c_2>0$ and $\alpha>d/2$ such that
	\[
	c_1\left(1+\lambda_\ell\right)^{-\alpha}
	\le
	C_\ell
	\le
	c_2\left(1+\lambda_\ell\right)^{-\alpha},
	\qquad
	\ell \ge 0,
	\]
	where $\lambda_\ell=\ell(\ell+d-1)$ are the eigenvalues of the
	Laplace--Beltrami operator on $\mathbb{S}^d$.
\end{condition}

Condition~\ref{cond:polydecay} characterizes Mat\'ern-type Gaussian field priors
on the sphere.
The requirement $\alpha>d/2$ ensures summability of the prior variances and
guarantees that the prior defines a centered Gaussian measure supported on
$L^2(\mathbb{S}^d)$.
Moreover, the decay rate of $\{C_\ell:\ell\ge0\}$ determines the Sobolev
regularity of prior sample paths and the reproducing kernel Hilbert space of the
Gaussian measure.

Polynomial decay of angular power spectra is also strongly motivated by
applications.
In particular, in the statistical analysis of the cosmic microwave background,
temperature and polarization fields are routinely modeled as isotropic Gaussian
fields with approximately polynomial spectral decay at high multipoles
\cite{Planck2018I,Planck2018V}.

This condition will be imposed throughout the contraction analysis and allows
for sharp bias--variance trade--offs to be derived explicitly.

\paragraph{Heuristic choice of the truncation level.}
Truncation at level \(L_n\) induces a squared bias of order
\[
\sum_{\ell>L_n}
\sum_{m=1}^{M_{d,\ell}} a_{0;\ell,m}^2
\simeq L_n^{-2\beta},
\qquad f_0\in H^\beta(\mathbb{S}^d),
\]
while retaining frequencies up to \(L_n\) yields an effective variance of order
\(L_n^d/n\), reflecting the growth of harmonic multiplicities on
\(\mathbb{S}^d\).
In the Bayesian setting, posterior variance is further shaped by the prior
through the spectral decay parameter \(\alpha\), which determines the strength
of regularization at high frequencies.
Balancing approximation bias and prior-regularized variance at a heuristically
leads to the choice
\[
L_n \simeq n^{1/(2\alpha+d)},
\]
which anticipates the rate later confirmed by the posterior contraction theorem
below.

\begin{remark}[Relation to the literature]
	Spectral truncation and shrinkage are central to the theory of nonparametric
	regression and inverse problems.
	In Euclidean settings, they underlie Pinsker estimators, smoothing splines,
	kernel ridge regression, and the classical Gaussian sequence model
	\cite{DonohoJohnstone94,Tsybakov09}.
	
	On the sphere and more general manifolds, analogous principles appear in
	needlet-based regression, thresholding estimators, and multiresolution methods
	\cite{BKMP09,Dur16,Mon11}.
	
	From a Bayesian perspective, truncated spectral posteriors can be interpreted
	as finite-dimensional approximations of Gaussian process priors.
	Explicit truncation renders these approximations transparent and ensures that
	the resulting estimators are both computationally feasible and theoretically
	well-aligned with classical bias--variance trade-offs
	\cite{GineNickl16,vanDerVaartvanZanten08}.
\end{remark}

\paragraph{Posterior contraction under truncation.}
We now state a posterior contraction result for the truncated Gaussian prior.
Exploiting the diagonal harmonic representation, the $L^2(\mathbb{S}^d)$ error
decomposes naturally into an approximation bias due to spectral truncation and
a variance term reflecting posterior uncertainty.
The interaction between these two components, governed by the growth of
Laplace--Beltrami eigenvalues and their multiplicities, determines both the
optimal truncation level and the resulting contraction rate.

\begin{theorem}[Posterior contraction under spectral truncation]
	\label{thm:contraction}
	Assume the nonparametric regression model of Section~\ref{sec:posterior} and let the prior
	satisfy Condition~\ref{cond:polydecay} for some $\alpha>d/2$.
	Assume that the true regression function satisfies
	\[
	f_0 \in H^\beta(\mathbb{S}^d),
	\qquad
	0<\beta\le\alpha.
	\]
	
	Let the truncation level be chosen as
	\[
	L_n \simeq n^{1/(2\alpha+d)}.
	\]
	Then the posterior distribution induced by the truncated Gaussian prior
	contracts around $f_0$ at rate
	\[
	\rho_n = n^{-\beta/(2\alpha+d)}
	\]
	in the $L^2(\mathbb{S}^d)$ norm, in the sense that for every sufficiently large
	constant $M>0$,
	\[
	\mathbb{P}_{f_0}^{(n)}\!\left(
	\Pi\!\left(
	f:\|f-f_0\|_{L^2(\mathbb{S}^d)}>M\rho_n
	\mid \bm{y},\bm{x}
	\right)
	\right)
	\longrightarrow 0
	\quad\text{as } n\to\infty.
	\]
\end{theorem}

Theorem~\ref{thm:contraction} establishes the rate at which the Bayesian
posterior concentrates around the true regression function under the assumed
prior and truncation scheme.
We now show that this rate is not only achievable but also information-theoretically optimal.

\begin{corollary}[Minimax optimality under correct prior calibration over Sobolev classes]
	\label{cor:minimax}
	Under the assumptions of Theorem~\ref{thm:contraction}, suppose in addition
	that the prior smoothness parameter satisfies $\alpha=\beta$.
	Then the posterior contraction rate
	\[
	\rho_n = n^{-\beta/(2\beta+d)}
	\]
	is minimax-optimal over Sobolev balls $H^\beta(\mathbb{S}^d)$ with respect to
	the $L^2(\mathbb{S}^d)$ loss.
	In particular, no estimator, Bayesian or frequentist, can achieve a uniformly
	faster rate over $H^\beta(\mathbb{S}^d)$.
\end{corollary}

\begin{remark}[Interpretation, calibration, and saturation]\label{rem:interpretation}
	The contraction rate in Theorem~\ref{thm:contraction} reflects the interaction
	between the Sobolev regularity $\beta$ of the true regression function and the
	regularization strength encoded by the prior smoothness parameter $\alpha$,
	through the decay of the angular power spectrum.
	Under the assumption $\beta\le\alpha$, the posterior contraction rate depends
	explicitly on the regularity of the truth and arises from a nontrivial balance
	between the approximation bias induced by spectral truncation and the posterior
	variance controlled by the prior.
	
	In the calibrated case $\alpha=\beta$, this balance is optimal and yields the
	minimax-optimal rate over the Sobolev class $H^\beta(\mathbb{S}^d)$.
	When $\beta<\alpha$, the prior is smoother than the truth and the contraction
	rate deteriorates accordingly, reflecting the increased bias induced by
	oversmoothing.
	
	Although not covered by Theorem~\ref{thm:contraction}, the same bias--variance
	decomposition shows that if the truth is smoother than the prior RKHS,
	$\beta>\alpha$, the posterior contraction rate saturates at the
	prior-limited rate
	\[
	n^{-\alpha/(2\alpha+d)},
	\]
	independently of the additional smoothness of $f_0$.
	In this regime, the procedure remains minimax-optimal over
	$H^\alpha(\mathbb{S}^d)$, but does not adapt to higher regularity of the truth.
	
	This saturation phenomenon is classical for fixed-order spectral Gaussian
	priors and smoothing splines, and explains why Theorem~\ref{thm:contraction}
	focuses on the regime $\beta\le\alpha$, where the dependence of the contraction
	rate on the unknown smoothness of the regression function is most transparent.
\end{remark}

\begin{remark}[Comparison with needlet regression and Besov regularity]
	The contraction rate in Theorem~\ref{thm:contraction} coincides with the optimal
	rates achieved by spectral cut-off and needlet-based estimators on
	$\mathbb{S}^d$, including thresholding and block-thresholding procedures
	\cite{BKMP09,Dur16,Mon11}.
	These estimators attain minimax-optimal rates over a broad range of Besov spaces
	on the sphere, which strictly contain Sobolev and H\"older classes.
	
	The present analysis is formulated in the spherical harmonic basis and is
	therefore naturally tailored to Sobolev regularity, corresponding to quadratic
	(energy-based) control of harmonic coefficients.
	In particular, H\"older classes on $\mathbb{S}^d$ embed continuously into
	Sobolev spaces of slightly lower smoothness, and the rate
	$\rho_n = n^{-\beta/(2\beta+d)}$ remains optimal for such classes.
	
	By contrast, adaptive estimation over Besov spaces typically requires localized
	multiscale representations, such as needlets, together with nonlinear
	thresholding.
	While the present Gaussian prior does not adapt automatically to unknown
	smoothness or sparsity, needlet-based constructions exploit spatial localization
	and scale-dependent shrinkage to achieve adaptation over a wide range of
	function classes.
	This parallels the Euclidean setting, where Fourier-based Gaussian priors yield
	sharp rates over Sobolev classes under correct calibration, whereas wavelet-based
	methods are required to attain adaptivity over Besov scales.
\end{remark}

\begin{remark}[Unknown prior smoothness and Whittle-type adaptation]
	\label{rem:unknown-alpha-whittle}
	Throughout the paper, the smoothness parameter $\alpha$ governing the decay of
	the prior angular power spectrum is assumed fixed and known, while the Sobolev
	regularity $\beta$ of the true regression function is unknown. This corresponds
	to the classical fixed-order Gaussian process regression setting and allows for a
	fully explicit analysis of posterior contraction rates, bias--variance trade-offs,
	and variational structure.
	
	Allowing the prior smoothness parameter $\alpha$ to be unknown and inferred from
	the data is a natural extension. From a Bayesian perspective, this may be
	approached either through hierarchical Gaussian process priors placing a
	hyperprior on $\alpha$, or through empirical Bayes procedures based on marginal
	likelihood maximization; see, for instance,
	\cite{vanDerVaartvanZanten08,vanDerVaartvanZanten09}. In the present spectral
	framework, the latter is closely related to Whittle-type likelihood methods and
	to semiparametric estimation of spectral regularity parameters for Gaussian
	random fields in the harmonic domain, as developed for spherical settings in
	\cite{dljber}.
	
	Bayesian counterparts of Whittle likelihood methods have been studied in the
	context of spectral density estimation for time series and Gaussian processes,
	see for example \cite{ChoGhoRoy04}. On the sphere, however, the presence
	of frequency-dependent harmonic multiplicities and the associated spectral
	geometry introduce additional analytical challenges, and a systematic treatment
	of hierarchical or Whittle-based adaptive Bayesian procedures on $\mathbb{S}^d$
	lies beyond the scope of the present work.
	
	We emphasize that the restriction to fixed $\alpha$ is not intrinsic to the
	proposed framework, but rather reflects a deliberate modeling choice aimed at
	obtaining exact harmonic diagonalization, closed-form posteriors, and sharp
	non-asymptotic guarantees. The development of fully adaptive spectral Bayesian
	methods on the sphere, and their precise connection with Whittle likelihood and
	semiparametric inference, is left for future research.
\end{remark}

\begin{remark}[Spectral coefficient representation and functional viewpoint]
	\label{rem:functional}
	
	The diagonal Gaussian sequence representation induced by spherical harmonics
	provides a natural bridge between Bayesian nonparametric regression on
	$\mathbb{S}^d$ and functional data analysis on manifolds.
	Rather than viewing the model purely in terms of pointwise noisy observations,
	one may interpret inference as acting directly on the sequence of harmonic
	coefficients.
	
	Under uniform random design, the regression model admits the representation \eqref{eqn:model};
	this is formally analogous to functional observation models studied in
	functional data analysis, where inference is performed coefficient-wise in an
	orthonormal basis of $L^2$; see, for example, \cite{Bosq00,PanaretosTavakoli13}.
	
	From this perspective, the Gaussian prior with angular power spectrum
	$\{C_\ell\}$ induces independent shrinkage on each harmonic block,
	\[
	a_{\ell,m} \sim \mathcal{N}(0,C_\ell),
	\qquad
	m=1,\dots,M_{d,\ell},
	\]
	while posterior contraction is governed by the aggregate contribution
	\[
	\sum_{\ell=0}^\infty
	\sum_{m=1}^{M_{d,\ell}}
	\mathbb{E}\!\left[
	\left(a_{\ell,m}-a_{0;\ell,m}\right)^2
	\mid \bm y,\bm x
	\right],
	\]
	which combines coefficient-wise estimation error with the geometric growth of
	harmonic multiplicities.
	
	Closely related ideas appear in recent work on spherical functional
	autoregressions, where stochastic dynamics are formulated directly in harmonic
	coordinates and asymptotic behavior is characterized through spectral norms;
	see \cite{cm21}.
	In those models, dependence across time acts diagonally in the spherical
	harmonic basis, while global error control is obtained by summing over
	frequency-dependent blocks, mirroring the structure of the present contraction
	analysis.
	
	More abstractly, both Bayesian regression and functional time series models on
	$\mathbb{S}^d$ can be viewed as acting on the infinite-dimensional coefficient
	vector
	\[
	\mathbf a
	=
	\left\{ a_{\ell,m} : \ell \ge 0,\; m=1,\dots,M_{d,\ell} \right\}
	\in \ell^2,
	\]
	equipped with a geometry-induced norm of the form
	\[
	\|\mathbf a\|_{\mathcal H_\alpha}^2
	=
	\sum_{\ell=0}^\infty
	(1+\lambda_\ell)^{\alpha}
	\sum_{m=1}^{M_{d,\ell}} a_{\ell,m}^2,
	\]
	which is equivalent to the Sobolev norm
	$\|f\|_{H^\alpha(\mathbb{S}^d)}^2$ under the harmonic isomorphism.
	
	In this representation, Bayesian posterior contraction corresponds to
	concentration of the posterior measure in $\mathcal H_\alpha$-balls around the
	true coefficient sequence, while functional autoregressive models on the sphere
	analyze stability and asymptotic behavior through analogous spectral norms.
	The common structure is diagonal inference at the coefficient level combined
	with global control through a Laplace--Beltrami–induced functional norm.
\end{remark}

\paragraph{Spherical Mat\'ern priors.}

A particularly important and widely used class of isotropic Gaussian field
priors on $\mathbb{S}^d$ is provided by spherical Mat\'ern fields.
These priors arise naturally through functional calculus of the
Laplace--Beltrami operator and admit a direct interpretation as solutions to
stochastic partial differential equations on the sphere.
As such, they provide a canonical link between spectral decay, Sobolev
regularity, and geometrically intrinsic smoothing.

\begin{definition}[Spherical Mat\'ern prior]
	Let $\alpha>d/2$ and $\kappa>0$.
	A spherical Mat\'ern prior is defined as the centered Gaussian field on
	$\mathbb{S}^d$ with spherical harmonic expansion
	\[
	f(x)
	=
	\sum_{\ell=0}^\infty
	\sum_{m=1}^{M_{d,\ell}}
	a_{\ell,m}\,Y_{\ell,m}(x),
	\qquad
	a_{\ell,m}
	\sim
	\mathcal{N}\!\left(
	0,
	\left(\kappa^2+\lambda_\ell\right)^{-\alpha}
	\right),
	\]
	where $\lambda_\ell=\ell(\ell+d-1)$ are the eigenvalues of the
	Laplace--Beltrami operator on $\mathbb{S}^d$.
\end{definition}

The corresponding angular power spectrum satisfies
\[
C_\ell \simeq \left(1+\lambda_\ell\right)^{-\alpha}.
\]
Thus spherical Mat\'ern priors satisfy Condition~\ref{cond:polydecay}.
In particular, typical prior sample paths belong almost surely to
$H^s(\mathbb{S}^d)$ for all $s<\alpha-d/2$, while the reproducing kernel Hilbert
space of the Gaussian measure is norm-equivalent to $H^\alpha(\mathbb{S}^d)$.

It is important to distinguish between the almost-sure regularity of prior sample
paths and the regularization strength induced by the prior.
The loss of $d/2$ derivatives between the RKHS smoothness $\alpha$ and the
sample path regularity $s<\alpha-d/2$ is a classical feature of Gaussian random
fields on $d$-dimensional manifolds.
Posterior contraction rates, however, are governed by the RKHS associated with
the prior rather than by the typical smoothness of prior draws.
Consequently, correct calibration of the prior corresponds to matching the RKHS
smoothness parameter $\alpha$ with the Sobolev regularity $\beta$ of the true
regression function, even though prior sample paths are almost surely rougher.

From a statistical perspective, spherical Mat\'ern priors play the same role on
$\mathbb{S}^d$ as classical Mat\'ern priors do in Euclidean Gaussian process
regression: they provide a one-parameter family of priors indexed by a
smoothness parameter $\alpha$, allowing regularity assumptions on the unknown
function to be encoded in a transparent and geometrically intrinsic way.

The scale parameter $\kappa$ controls the correlation length and low-frequency
behavior of the prior but does not affect posterior contraction rates, which are
governed by the high-frequency spectral decay encoded by $\alpha$.
Accordingly, $\kappa$ enters only through multiplicative constants and is
suppressed in rate statements.

\begin{corollary}[Posterior contraction for spherical Mat\'ern priors]
	\label{cor:matern}
	Let the prior be a spherical Mat\'ern prior with smoothness parameter
	$\alpha>d/2$, and assume that the true regression function satisfies
	$f_0\in H^\beta(\mathbb{S}^d)$ for some $0<\beta\le\alpha$.
	Then, with truncation level
	\[
	L_n \simeq n^{1/(2\alpha+d)},
	\]
	the posterior distribution contracts around $f_0$ at rate
	\[
	\rho_n = n^{-\beta/(2\alpha+d)}
	\]
	in the $L^2(\mathbb{S}^d)$ norm.
\end{corollary}

The present Bayesian contraction result mirrors the spectral efficiency
properties of Whittle likelihood--based estimators for Mat\'ern-type Gaussian
fields on compact manifolds, while extending these ideas to Bayesian
nonparametric regression and posterior contraction.
In particular, Mat\'ern-type priors on the sphere yield minimax-optimal rates
and admit a transparent spectral interpretation closely aligned with classical
results in Euclidean Gaussian process theory
\cite{KorteStapffKarvonenMoulines25}.

\begin{remark}[Fixed-order nature of spherical Mat\'ern priors]
	The behavior of spherical Mat\'ern priors described in
	Corollary~\ref{cor:matern} is a direct specialization of the general
	bias--variance and calibration principles discussed in
	Remark~\ref{rem:interpretation}.
	In particular, spherical Mat\'ern priors correspond to polynomial spectral
	priors with
	$C_\ell \simeq (1+\lambda_\ell)^{-\alpha}$ and therefore act as fixed-order
	Sobolev regularizers on $\mathbb{S}^d$.
	
	As in Euclidean Gaussian process regression with Mat\'ern kernels, such priors
	achieve minimax-optimal rates under correct calibration of the smoothness
	parameter, but do not adapt automatically to unknown regularity.
	This fixed-order behavior is classical for Mat\'ern Gaussian process priors
	and smoothing splines in Euclidean settings
	\cite{GineNickl16,vanDerVaartvanZanten08}, and has close analogues in spectral
	and Whittle-type inference for isotropic Gaussian fields on manifolds
	\cite{dljber,KorteStapffKarvonenMoulines25}.	
\end{remark}

\begin{remark}[Spectral regularization, kernels, and Whittle likelihood]
	\label{rem:spectral-whittle}
	The explicit harmonic characterization of the posterior mean reveals that
	Bayesian regression on $\mathbb{S}^d$ is fundamentally a form of spectral
	regularization driven by the Laplace--Beltrami operator.
	In harmonic coordinates, posterior inference amounts to frequency-dependent
	shrinkage of empirical coefficients, with penalties determined by the angular
	power spectrum $\{C_\ell\}$.
	
	This mechanism admits several equivalent interpretations.
	First, the posterior mean coincides with a ridge-type estimator in the
	spherical harmonic basis, with penalty weights proportional to
	$\sigma^2/(nC_\ell)$.
	For polynomially decaying spectra
	$C_\ell \simeq (1+\lambda_\ell)^{-\alpha}$, this penalty is equivalent to a
	Sobolev norm induced by fractional powers of the Laplace--Beltrami operator,
	yielding a geometrically intrinsic smoothing spline on the sphere
	\cite{BelkinNiyogiSindhwani06,wahba90}.
	
	Equivalently, the isotropic Gaussian prior induces a rotation-invariant
	reproducing kernel, and the posterior mean coincides with kernel ridge
	regression on $\mathbb{S}^d$.
	For Mat\'ern-type spectra, the associated RKHS norm is norm-equivalent to a
	Sobolev space on the sphere, showing that Bayesian regression, smoothing
	splines, and kernel methods are mathematically identical procedures expressed
	in different languages.
	
	The same Laplace--Beltrami spectral structure also underlies Whittle
	likelihood--based methods for isotropic Gaussian random fields on the sphere.
	In such approaches, the likelihood is approximated in the eigenbasis of the
	Laplace--Beltrami operator, leading to diagonal or asymptotically diagonal
	likelihoods in spherical harmonic or needlet coordinates.
	Needlet Whittle estimators \cite{dlmejs} and Gaussian semiparametric spectral
	methods \cite{dljber} exploit the high-frequency behavior of $C_\ell$ to
	estimate smoothness and scale parameters, while recent work on compact
	Riemannian manifolds \cite{KorteStapffKarvonenMoulines25} shows that Mat\'ern-type
	spectral penalties yield semiparametrically optimal inference for smoothness
	parameters.
	
	While Whittle and semiparametric approaches focus on likelihood-based
	inference for covariance and smoothness parameters, the present work is
	concerned with posterior contraction and uncertainty quantification for
	function estimation.
	Both frameworks, however, rely on the same intrinsic spectral regularization
	principle induced by the Laplace--Beltrami operator.
\end{remark}

\section{Variational Interpretation}
\label{sec:variational}

This section provides a variational interpretation of Bayesian regression on
$\mathbb{S}^d$.
Rather than viewing Gaussian process regression purely as a probabilistic
procedure, we show that posterior inference induces a geometrically intrinsic
regularization mechanism governed by the Laplace--Beltrami operator.
This perspective clarifies the role of the angular power spectrum, establishes
an exact equivalence with penalized least-squares estimation on the sphere, and
makes explicit the variational structure that underlies the Bayesian regression
procedures analyzed above.

Throughout this section, we work with the truncated harmonic space
$\mathcal{H}_{L_n}$ introduced in Section~\ref{sec:posterior}, and we denote by
$\hat f_{n,L_n}$ the posterior mean under the truncated Gaussian prior.


A central feature of Gaussian process regression in Euclidean settings is the
equivalence between posterior means and solutions of penalized least-squares
problems, a fact that underlies the theory of smoothing splines and kernel ridge
regression.
On $\mathbb{S}^d$, however, this equivalence is far from automatic: it must
respect the intrinsic geometry of the domain, the spectral structure of the
Laplace--Beltrami operator, and the growing multiplicities of spherical harmonic
eigenspaces.

In this subsection we show that, under uniform random design and isotropic
Gaussian field priors, Bayesian regression on the sphere admits an exact and
fully intrinsic variational characterization.
The resulting estimator is a geometrically natural ridge-type regularization
scheme, with penalty weights determined explicitly by the angular power
spectrum of the prior.
This provides a deterministic interpretation of posterior inference and
connects the probabilistic analysis of previous sections to classical
regularization theory on manifolds.

\begin{theorem}[Variational characterization of the posterior mean]
	\label{thm:variational}
	Consider the nonparametric regression model on $\mathbb{S}^d$ under uniform
	random design with noise variance $\sigma^2$.
	Let $\Pi$ be a centered isotropic Gaussian field prior with angular power
	spectrum $\{C_\ell:\ell\ge0\}$, truncated at level $L_n$, and let
	$\hat f_{n,L_n}$ denote the posterior mean under the resulting truncated
	posterior distribution.
	
	Then $\hat f_{n,L_n}$ is the unique minimizer over the finite-dimensional
	harmonic space $\mathcal{H}_{L_n}$ of the functional
	\[
	\mathcal{J}_n(f)
	=
	\frac{1}{n}\sum_{i=1}^n \bigl(y_i - f(x_i)\bigr)^2
	+
	\frac{\sigma^2}{n}
	\sum_{\ell=0}^{L_n}
	\sum_{m=1}^{M_{d,\ell}}
	\frac{a_{\ell,m}^2}{C_\ell},
	\]
	where $\{a_{\ell,m}:\ell \geq 0; m=1,\ldots,M_{d,\ell}\}$ are the spherical harmonic coefficients of $f$.
	
	Equivalently, $\hat f_{n,L_n}$ solves
	\[
	\hat f_{n,L_n}
	=
	\arg\min_{f\in\mathcal{H}_{L_n}} \mathcal{J}_n(f),
	\]
	and this minimizer exists and is unique.
\end{theorem}

Theorem~\ref{thm:variational} shows that Bayesian regression on $\mathbb{S}^d$
induces an intrinsic quadratic regularization mechanism governed by the
Laplace--Beltrami spectrum.
The penalty term acts diagonally in harmonic coordinates and suppresses
high-frequency components according to the prior spectrum $\{C_\ell\}$, while
automatically accounting for the geometric multiplicities
$M_{d,\ell}$ of the eigenspaces.
In particular, for polynomially decaying spectra
$C_\ell \simeq (1+\lambda_\ell)^{-\alpha}$, the penalty is equivalent to a
Sobolev-type norm defined through fractional powers of the Laplace--Beltrami
operator.

This variational characterization is exact, non-asymptotic, and does not rely
on kernel representations or continuum limits.
It provides the spherical analogue of classical results for smoothing splines
and kernel ridge regression in Euclidean domains, while making explicit the
distinctive role of geometry through the harmonic structure of $\mathbb{S}^d$. 
The proof, which exploits the diagonal Gaussian sequence representation derived
earlier and the strict convexity of $\mathcal{J}_n$, is given in
Section~\ref{sec:proofs}.

\begin{remark}[Variational origin of posterior contraction]
	The contraction rate in Theorem~\ref{thm:contraction} can be interpreted
	directly through the variational formulation of the posterior mean.
	By Theorem~\ref{thm:variational}, $\hat f_{n,L_n}$ minimizes a penalized
	least-squares functional whose penalty suppresses high-frequency spherical
	harmonics according to the angular power spectrum $\{C_\ell\}$.
	
	The proof of Theorem~\ref{thm:contraction} exploits this structure
	by decomposing the posterior risk into a truncation bias term, controlled by
	the Sobolev regularity of $f_0$, and a variance term arising from stochastic
	fluctuations of the empirical coefficients.
	The optimal choice $L_n\simeq n^{1/(2\alpha+d)}$ corresponds precisely to the
	point at which the deterministic regularization induced by the penalty
	matches the noise level in the Gaussian sequence model.
	
	In this sense, posterior contraction on $\mathbb{S}^d$ is a direct
	consequence of the spectral regularization mechanism encoded in the
	variational formulation.
\end{remark}

The variational characterization in Theorem~\ref{thm:variational} holds for
arbitrary isotropic Gaussian priors through their angular power spectrum
$\{C_\ell\}$.
A particularly important specialization is obtained when the spectrum is
generated by functional calculus of the Laplace--Beltrami operator.
In this case, the abstract frequency-dependent penalty admits a direct
geometric interpretation in terms of differential operators on
$\mathbb{S}^d$, leading to an intrinsic notion of smoothing.
This situation arises precisely for spherical Mat\'ern priors, which play a
distinguished role both in Gaussian process regression and in spectral
likelihood-based inference on manifolds.

\begin{corollary}[Variational form of spherical Mat\'ern smoothing]
	\label{cor:variational-matern}
	Let the prior be a spherical Mat\'ern field with parameters
	$\alpha>d/2$ and $\kappa>0$, so that
	\[
	C_\ell=(\kappa^2+\lambda_\ell)^{-\alpha},
	\qquad
	\lambda_\ell=\ell(\ell+d-1).
	\]
	Then the posterior mean $\hat f_{n,L_n}$ is the unique minimizer over
	$\mathcal{H}_{L_n}$ of
	\[
	\frac{1}{n}\sum_{i=1}^n (y_i-f(x_i))^2
	+
	\frac{\sigma^2}{n}
	\sum_{\ell=0}^{L_n}\sum_{m=1}^{M_{d,\ell}}
	(\kappa^2+\lambda_\ell)^{\alpha} a_{\ell,m}^2.
	\]
	
	Equivalently, $\hat f_{n,L_n}$ is a Laplace--Beltrami smoothing spline of
	order $\alpha$ on $\mathbb{S}^d$.
\end{corollary}

\begin{remark}[Non-adaptivity and smoothing behavior of Mat\'ern priors]
	Spherical Mat\'ern priors act as fixed-order Sobolev regularizers determined
	by the smoothness parameter $\alpha$.
	As a consequence, they achieve minimax-optimal rates when the regularity of
	the true regression function satisfies $\beta\le\alpha$, but do not adapt
	automatically to unknown smoothness.
	
	If $\alpha>\beta$, the estimator oversmooths: high-frequency components of
	the signal are excessively penalized, and the contraction rate is driven by
	the bias term.
	If $\alpha<\beta$, the estimator undersmooths: variance dominates and the
	full regularity of the signal is not exploited.
	
	This behavior mirrors the classical bias--variance trade-off for Mat\'ern
	priors in Euclidean Gaussian process regression and for Whittle likelihood
	estimation of smoothness parameters on manifolds.
	In both settings, Mat\'ern models provide sharp inference when correctly
	calibrated, but lack adaptivity without additional hierarchical or
	multiscale structure.
	
	From a variational viewpoint, this reflects the fact that the penalty operator
	$(\kappa^2-\Delta_{\mathbb{S}^d})^{\alpha}$ enforces a single Sobolev scale.
	Adaptive estimation over Besov or spatially inhomogeneous classes typically
	requires localized bases, such as needlets, together with nonlinear or
	hierarchical regularization.
\end{remark}

\begin{remark}[Intrinsic geometric regularization]
	\label{rem:geometric_regularization}
	The variational penalty induced by a spherical Mat\'ern prior admits a direct
	interpretation in terms of Sobolev regularization on $\mathbb{S}^d$.\\
	When the angular power spectrum exhibits polynomial decay,
	\[
	C_\ell \simeq (1+\lambda_\ell)^{-\alpha},
	\qquad \alpha > 0,
	\]
	the penalty term in the variational functional $\mathcal{J}_n$ admits a direct
	intrinsic geometric interpretation.
	Indeed, for any $f\in\mathcal{H}_{L_n}$ with spherical harmonic coefficients
	$\{a_{\ell,m}\}$,
	\[
	\sum_{\ell=0}^{L_n}\sum_{m=1}^{M_{d,\ell}}
	\frac{a_{\ell,m}^2}{C_\ell}
	\;\simeq\;
	\sum_{\ell=0}^{L_n}\sum_{m=1}^{M_{d,\ell}}
	(1+\lambda_\ell)^{\alpha} a_{\ell,m}^2
	\;=\;
	\|f\|_{H^{\alpha}(\mathbb{S}^d)}^2,
	\]
	where $\lambda_\ell=\ell(\ell+d-1)$ are the eigenvalues of the Laplace--Beltrami
	operator on $\mathbb{S}^d$.
	
	Consequently, the posterior mean minimizes an empirical risk penalized by a
	Sobolev norm defined intrinsically on the sphere.
	Equivalently, Gaussian process regression on $\mathbb{S}^d$ induces smoothing
	through fractional powers of the Laplace--Beltrami operator.
	This regularization mechanism is rotation invariant, geometrically intrinsic,
	and independent of any embedding of the sphere into Euclidean space.
	
	From a variational perspective, the estimator $\hat f_{n,L_n}$ can therefore be
	interpreted as a spherical smoothing spline, extending the classical theory of
	Sobolev penalization and spline smoothing to compact manifolds.
	At the same time, the associated Gaussian prior induces a reproducing kernel
	Hilbert space whose norm is equivalent to $H^{\alpha}(\mathbb{S}^d)$, showing
	that Bayesian regression on the sphere coincides exactly with Laplace--Beltrami
	kernel ridge regression with a Mat\'ern-type kernel.
	
	This interpretation also clarifies the connection with spectral and Whittle-type
	methods for isotropic Gaussian random fields on the sphere.
	In those approaches, regularization is likewise governed by powers of
	$(\kappa^2+\lambda_\ell)$ and by the high-frequency behavior of the angular
	power spectrum, although the inferential target is covariance structure rather
	than function recovery.
\end{remark}

\begin{remark}[Connection with smoothing splines and  kernel ridge regression]
	\label{rem:splines_krr}
	
	The variational characterization of the posterior mean places spherical Gaussian
	process regression within the classical theory of smoothing splines and kernel
	ridge regression, formulated intrinsically on $\mathbb{S}^d$.
	
	The isotropic Gaussian prior induces the reproducing kernel
	\[
	K(x,x')
	=
	\sum_{\ell=0}^\infty
	C_\ell \frac{M_{d,\ell}}{\omega_d}
	G_\ell^{((d-1)/2)}(\langle x,x' \rangle),
	\]
	which is rotation invariant, and 
	where $G_\ell^{((d-1)/2)}$ are Gegenbauer polynomials and $\omega_d$ denotes the
	surface area of $\mathbb{S}^d$.
	Thus isotropic kernels on the sphere depend only on geodesic distance through
	$\langle x,x' \rangle$ and are completely characterized by the angular power
	spectrum.
	
	Restricting $K$ to the truncated space $\mathcal{H}_{L_n}$ yields a finite-rank
	kernel $K_{L_n}$, and the posterior mean $\hat f_{n,L_n}$ coincides with the kernel
	ridge regression estimator
	\[
	\hat f_{n,L_n}
	=
	\arg\min_{f \in \mathcal{H}_{L_n}}
	\left\{
	\frac{1}{n} \sum_{i=1}^n (y_i - f(x_i))^2
	+
	\lambda_n \|f\|_{\mathcal{H}_K}^2
	\right\},
	\qquad
	\lambda_n = \sigma^2 / n.
	\]
	
	For Mat\'ern-type spectra, the RKHS norm $\|\cdot\|_{\mathcal{H}_K}$ is equivalent
	to the Sobolev norm $H^\alpha(\mathbb{S}^d)$.
	Hence Bayesian regression, kernel ridge regression, and smoothing spline
	estimators on the sphere are the same mathematical object expressed in different
	languages.
	
	This equivalence is exact, non-asymptotic, and geometrically intrinsic.
	It clarifies the role of the Laplace--Beltrami spectrum and its multiplicities in
	regularization and bias--variance trade-offs, and explains why Gaussian process
	priors act as spectral regularizers on compact manifolds.
\end{remark}

\begin{remark}[Variational structure and functional spectral interpretation]
	\label{rem:functvar}
	
	The variational characterization of Remark ~\ref{rem:functional} provides a
	functional counterpart to the coefficient-wise Bayesian analysis of
	Section~\ref{sec:posterior}.
	The posterior mean $\hat f_{n,L_n}$ is the unique minimizer over
	$\mathcal H_{L_n}$ of the penalized functional
	\[
	\frac{1}{n}\sum_{i=1}^n \left(y_i - f(x_i)\right)^2
	+
	\frac{\sigma^2}{n}
	\left\|
	\psi\!\left(-\Delta_{\mathbb{S}^d}\right)^{1/2} f
	\right\|_{L^2(\mathbb{S}^d)}^2,
	\]
	where the penalty is defined through spectral calculus of the
	Laplace--Beltrami operator.
	
	This shows that Bayesian regression on $\mathbb{S}^d$ combines diagonal
	coefficient-wise inference in harmonic coordinates with a single intrinsic
	functional norm controlling global smoothness.
	The regularization acts through the geometry of the sphere rather than through
	coordinate-wise constraints.
	
	This perspective suggests natural extensions to genuinely functional spherical
	models, such as time-indexed curve- or field-valued observations, where
	dependence may be modeled spectrally while regularity is enforced through
	operator-based norms.
	A systematic Bayesian treatment of such models is left for future work.
\end{remark}

\section{Numerical experiments}
\label{sec:num}

This section provides numerical evidence supporting the posterior contraction
results established in Section~\ref{sec:posterior}.
We consider Bayesian nonparametric regression on the unit sphere
$\mathbb{S}^2$ under uniform random design, using truncated spherical
Mat\'ern priors and evaluating performance in the $L^2(\mathbb{S}^2)$ norm.

\subsection{Experimental setup}

The data are generated according to the regression model
\[
y_i = f_0(x_i) + \varepsilon_i,
\qquad
\varepsilon_i \sim \mathcal{N}(0,\sigma^2),
\]
with design points $x_i$ sampled independently and uniformly from
$\mathbb{S}^2$.
The true regression function $f_0$ is constructed as a finite spherical harmonic
expansion up to degree $L_0=10$, with coefficients decaying as
$(1+\lambda_\ell)^{-\beta/2}$ for $\beta=2$, ensuring that
$f_0 \in H^\beta(\mathbb{S}^2)$.

We place a truncated isotropic Gaussian prior on $f$ with angular power spectrum
\[
C_\ell = (1+\lambda_\ell)^{-\alpha},
\qquad \alpha=2,
\]
corresponding to a spherical Mat\'ern prior.
The truncation level is chosen as
\[
L_n = \left\lfloor c\, n^{1/(2\alpha+2)} \right\rfloor,
\qquad c=2.5,
\]
in accordance with the theoretical bias--variance balance.
For each sample size $n$, the posterior mean is computed explicitly in harmonic
coordinates and evaluated on a fine spherical grid.
Each experiment is repeated independently $50$ times.

The root mean squared error (RMSE) is computed as
\[
\mathrm{RMSE}
=
\left(
\int_{\mathbb{S}^2}
\bigl(\hat f_{n,L_n}(x) - f_0(x)\bigr)^2 \, dx
\right)^{1/2},
\]
approximated numerically on the grid and averaged over repetitions.

\subsection{Empirical contraction rates}

Table~\ref{tab:rmse} reports the empirical RMSE together with the corresponding
truncation levels $L_n$ for increasing sample sizes.
A clear monotone decay of the RMSE is observed as $n$ increases, while the
truncation level grows slowly, reflecting the polynomial growth predicted by
Theorem~\ref{thm:contraction}.

\begin{table}[h]
	\centering
	\begin{tabular}{ccc}
		\hline
		$n$ & $L_n$ & RMSE \\
		\hline
		50   & 4 & $2.45\times10^{-1}$ \\
		100  & 5 & $2.04\times10^{-1}$ \\
		200  & 6 & $1.74\times10^{-1}$ \\
		400  & 6 & $1.27\times10^{-1}$ \\
		800  & 7 & $1.03\times10^{-1}$ \\
		1600 & 8 & $8.41\times10^{-2}$ \\
		3200 & 9 & $6.90\times10^{-2}$ \\
		\hline
	\end{tabular}
	\caption{Empirical $L^2(\mathbb{S}^2)$ RMSE of the posterior mean and corresponding
		spectral truncation level $L_n$ for increasing sample size $n$.}
	\label{tab:rmse}
\end{table}

\begin{figure}[t]
	\centering
	\includegraphics[width=\textwidth]{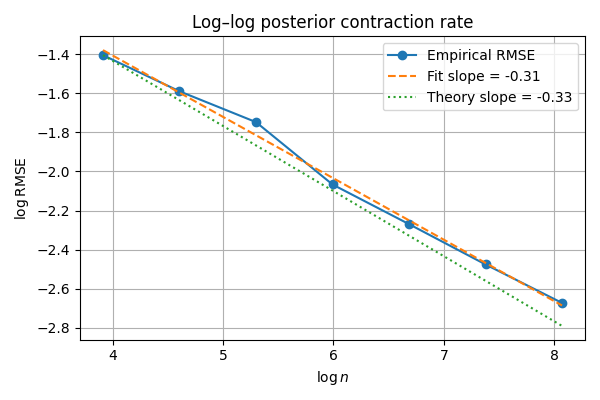}
	\caption{Log--log plot of the empirical RMSE of the posterior mean versus sample
		size $n$.
		The dashed line corresponds to the least-squares fit, while the dotted line
		indicates the theoretical slope $-\beta/(2\alpha+2)=-1/3$.}
	\label{fig:loglog}
\end{figure}

Figure~\ref{fig:loglog} displays the log--log plot of RMSE versus sample size.
A linear regression yields an estimated slope of $-0.314$, in close agreement
with the theoretical contraction rate
\[
-\frac{\beta}{2\alpha+2} = -\frac{1}{3}.
\]
The slight discrepancy is attributable to finite-sample effects and the use of a
fixed truncation constant.
Overall, the numerical results are consistent with the rate
predicted by the theory.

\subsection{Spatial structure of the posterior mean}

\begin{figure}[t]
	\centering
	\includegraphics[width=0.95\textwidth]{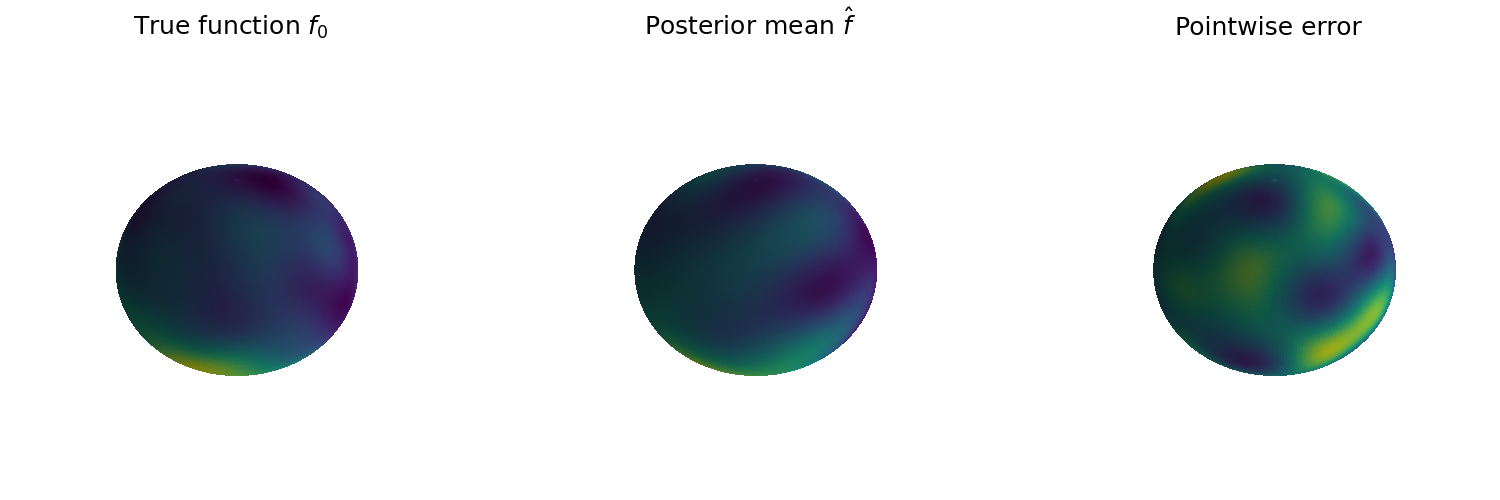}
	\caption{Visualization on $\mathbb{S}^2$ of the true regression function $f_0$
		(left), the posterior mean $\hat f_{n,L_n}$ (center), and the pointwise error
		$\hat f_{n,L_n}-f_0$ (right), for a representative realization.
		Color intensity represents function values on the sphere.}
	\label{fig:sphere}
\end{figure}

Figure~\ref{fig:sphere} shows the true regression function $f_0$, the posterior
mean $\hat f_{n,L_n}$, and the pointwise error
$\hat f_{n,L_n}-f_0$, visualized on $\mathbb{S}^2$ for a representative sample
size.
The posterior mean accurately recovers the large-scale structure of the signal,
while attenuating high-frequency oscillations in accordance with the spectral
regularization induced by the Mat\'ern prior.
The error surface exhibits no systematic spatial bias and is dominated by
localized fluctuations.

The numerical experiments thus corroborate the theoretical analysis:
Bayesian regression on the sphere with polynomially decaying angular power
spectra yields estimators that achieve the predicted contraction rates and
exhibit geometrically intrinsic smoothing behavior.

\subsection{Effect of prior mis-calibration}
\label{subsec:miscalibration}

We conclude the numerical study by investigating the effect of prior
mis-calibration, that is, the impact of choosing a smoothness parameter
$\alpha$ that does not coincide with the true Sobolev regularity
$\beta$ of the regression function.
This experiment illustrates concretely the non-adaptive nature of
polynomially decaying Gaussian priors and complements the theoretical
discussion in Sections~\ref{sec:posterior} and~\ref{sec:variational}.

We fix the true smoothness at $\beta=2$ and consider three prior choices:
$\alpha=1$ (undersmoothing), $\alpha=2$ (correct calibration), and
$\alpha=3$ (oversmoothing).
For each value of $\alpha$, the truncation level is chosen as
$L_n \simeq n^{1/(2\alpha+2)}$, in accordance with the theoretical bias--variance
balance.
Table~\ref{tab:miscalibration} reports the empirical $L^2(\mathbb{S}^2)$ RMSE of
the posterior mean for increasing sample sizes, together with the corresponding
truncation levels.

\begin{table}[t]
	\centering
	\begin{tabular}{c|cc|cc|cc}
		\hline
		& \multicolumn{6}{c}{Prior smoothness $\alpha$ (truth $\beta=2$)} \\
		$n$
		& \multicolumn{2}{c}{$\alpha=1$}
		& \multicolumn{2}{c}{$\alpha=2$}
		& \multicolumn{2}{c}{$\alpha=3$} \\
		& $L_n$ & RMSE & $L_n$ & RMSE & $L_n$ & RMSE \\
		\hline
		50   & 6  & $3.87\times10^{-1}$ & 4 & $2.45\times10^{-1}$ & 3 & $1.69\times10^{-1}$ \\
		100  & 7  & $3.13\times10^{-1}$ & 4 & $1.80\times10^{-1}$ & 4 & $1.41\times10^{-1}$ \\
		200  & 8  & $2.51\times10^{-1}$ & 5 & $1.56\times10^{-1}$ & 4 & $1.18\times10^{-1}$ \\
		400  & 10 & $2.13\times10^{-1}$ & 6 & $1.27\times10^{-1}$ & 4 & $9.66\times10^{-2}$ \\
		800  & 12 & $1.77\times10^{-1}$ & 7 & $1.08\times10^{-1}$ & 5 & $7.95\times10^{-2}$ \\
		1600 & 14 & $1.47\times10^{-1}$ & 7 & $8.10\times10^{-2}$ & 5 & $6.70\times10^{-2}$ \\
		3200 & 17 & $1.27\times10^{-1}$ & 8 & $6.45\times10^{-2}$ & 6 & $5.10\times10^{-2}$ \\
		\hline
		Theoretical slope
		& \multicolumn{2}{c}{$-1/2$}
		& \multicolumn{2}{c}{$-1/3$}
		& \multicolumn{2}{c}{$-1/4$} \\
		\hline
	\end{tabular}
	\caption{Empirical $L^2(\mathbb{S}^2)$ RMSE of the posterior mean under prior
		mis-calibration for true smoothness $\beta=2$.
		Theoretical slopes correspond to $-\beta/(2\alpha+2)$.}
	\label{tab:miscalibration}
\end{table}

\begin{figure}[t]
	\centering
	\includegraphics[width=0.95\textwidth]{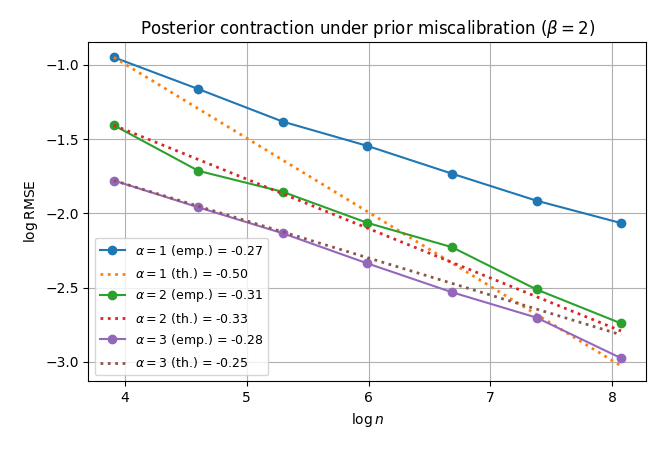}
	\caption{Log--log plot of the empirical RMSE versus sample size under prior
		mis-calibration.
		Solid lines correspond to empirical RMSE curves, while dotted lines indicate
		theoretical slopes $-\beta/(2\alpha+2)$ for $\alpha=1,2,3$.
		The correctly calibrated case $\alpha=\beta=2$ achieves the optimal rate.}
	\label{fig:miscalibration}
\end{figure}

Figure~\ref{fig:miscalibration} displays the corresponding log--log plots of the
RMSE as a function of the sample size.
In each case, the empirical slope closely matches the theoretical prediction
$-\beta/(2\alpha+2)$, confirming the contraction rates derived in
Theorem~\ref{thm:contraction}.
The correctly calibrated prior $\alpha=\beta$ yields the fastest admissible
rate $n^{-1/3}$, while both under- and over-smoothing lead to slower convergence.

These results provide clear numerical evidence of the non-adaptive behavior of
fixed-order polynomial spectral priors.
While spherical Mat\'ern priors yield sharp minimax rates when correctly
calibrated, misspecification of the smoothness parameter $\alpha$ leads to a
systematic degradation of performance, exactly as predicted by the theory.
This mirrors classical phenomena observed for Mat\'ern Gaussian process
regression in Euclidean settings and highlights the intrinsic bias--variance
trade-off induced by Laplace--Beltrami spectral regularization on the sphere.

\section{Proofs}
\label{sec:proofs}

\begin{proof}[Proof of Theorem~\ref{thm:contraction}]
We establish posterior contraction in the $L^2(\mathbb{S}^d)$--norm by exploiting
the exact diagonal Gaussian sequence representation induced by spherical harmonic
coordinates under uniform random design.
This representation allows the posterior risk to be decomposed explicitly into
a bias component arising from spectral truncation and prior-induced shrinkage,
and a variance component governed by the angular power spectrum and the
multiplicities of Laplace--Beltrami eigenspaces.
The proof proceeds by sharp control of these two contributions and by identifying
the truncation level at which they balance, taking full account of the geometric
growth of harmonic dimensions on $\mathbb{S}^d$.	

\medskip
\textit{Step 1: Orthogonal decomposition of the $L^2$--error.}
Since $\{Y_{\ell,m}\}_{\ell\ge0,\;1\le m\le M_{d,\ell}}$ forms an orthonormal basis
of $L^2(\mathbb{S}^d)$, the difference $f-f_0$ admits the expansion
\[
f(x)-f_0(x)
=
\sum_{\ell=0}^{\infty}
\sum_{m=1}^{M_{d,\ell}}
\left(
a_{\ell,m}-a_{0;\ell,m}
\right)
Y_{\ell,m}(x),
\]
with convergence in $L^2(\mathbb{S}^d)$.
By Parseval’s identity, it follows that
\[
\| f - f_0 \|_{L^2(\mathbb{S}^d)}^2
=
\sum_{\ell=0}^{\infty}
\sum_{m=1}^{M_{d,\ell}}
\left(
a_{\ell,m} - a_{0;\ell,m}
\right)^2.
\]
Under the truncated prior, we impose $a_{\ell,m}=0$ for $\ell>L_n$.
Accordingly, the $L^2$--error decomposes as
\[
\| f - f_0 \|_{L^2(\mathbb{S}^d)}^2
=
\sum_{\ell=0}^{L_n}
\sum_{m=1}^{M_{d,\ell}}
\left(
a_{\ell,m} - a_{0;\ell,m}
\right)^2
+
\sum_{\ell>L_n}
\sum_{m=1}^{M_{d,\ell}}
a_{0;\ell,m}^2,
\]
where the second term corresponds to the truncation bias.
For $0 \le \ell \le L_n$ the prior assigns independent Gaussian distributions
\[
a_{\ell,m} \sim \mathcal{N}(0,C_\ell).
\]
Combined with the Gaussian sequence representation of the likelihood in harmonic
coordinates and Gaussian conjugacy, this implies that the posterior distribution
factorizes over $(\ell,m)$ and satisfies
\[
a_{\ell,m} \mid \bm{y},\bm{x}
\sim
\mathcal{N}\left(
\mu_{\ell,m}^{(n)},\, v_\ell^{(n)}
\right),
\qquad
0 \le \ell \le L_n,
\]
with
\[
\mu_{\ell,m}^{(n)}
=
\frac{n C_\ell}{n C_\ell + \sigma^2}\,\hat a_{\ell,m},
\qquad
v_\ell^{(n)}
=
\frac{C_\ell \sigma^2}{n C_\ell + \sigma^2}.
\]

	\medskip
	\textit{Step 2: Posterior mean squared error decomposition.}
	We now take posterior expectation of the $L^2(\mathbb{S}^d)$--error conditionally
	on the observed data $(\bm{y},\bm{x})$.
	Using the orthogonal decomposition established in Step~1 and linearity of
	expectation, we obtain
	\[
	\mathbb{E}\!\left[
	\| f - f_0 \|_{L^2(\mathbb{S}^d)}^2
	\mid \bm{y},\bm{x}
	\right]
	=
	\sum_{\ell=0}^{L_n}
	\sum_{m=1}^{M_{d,\ell}}
	\mathbb{E}\!\left[
	\left( a_{\ell,m} - a_{0;\ell,m} \right)^2
	\mid \bm{y},\bm{x}
	\right]
	+
	\sum_{\ell>L_n}
	\sum_{m=1}^{M_{d,\ell}}
	a_{0;\ell,m}^2.
	\]
	The second term is deterministic and corresponds to the truncation bias, since
	$a_{\ell,m}=0$ almost surely under the truncated prior for $\ell>L_n$.
	
	For $0 \le \ell \le L_n$ and $m=1,\ldots, M_{d,\ell}$, the posterior distribution of each coefficient
	$a_{\ell,m}$ is Gaussian with mean $\mu_{\ell,m}^{(n)}$ and variance
	$v_\ell^{(n)}$.
	Moreover, the posterior factorizes over $(\ell,m)$.
For each retained harmonic mode, conditionally on $(\bm y,\bm x)$, we write
\[
a_{\ell,m} = \mu_{\ell,m}^{(n)} + Z_{\ell,m},
\]
where $Z_{\ell,m}$ is a centered Gaussian random variable satisfying
\[
\mathbb E[Z_{\ell,m}\mid \bm y,\bm x]=0,
\qquad
\mathbb E[Z_{\ell,m}^2\mid \bm y,\bm x]=v_\ell^{(n)}.
\]
Therefore,
\[
a_{\ell,m} - a_{0;\ell,m}
=
\mu_{\ell,m}^{(n)} - a_{0;\ell,m} + Z_{\ell,m}.
\]
Expanding the square and taking conditional expectation yields
\[
\begin{aligned}
	\mathbb{E}\!\left[
	( a_{\ell,m} - a_{0;\ell,m} )^2
	\mid \bm y,\bm x
	\right]
	&=
	( \mu_{\ell,m}^{(n)} - a_{0;\ell,m} )^2
	+ 2 ( \mu_{\ell,m}^{(n)} - a_{0;\ell,m} ) \mathbb E[Z_{\ell,m}] 
	+ \mathbb E[Z_{\ell,m}^2] \\
	&=
	( \mu_{\ell,m}^{(n)} - a_{0;\ell,m} )^2
	+ v_\ell^{(n)},
\end{aligned}
\]
since $\mathbb E[Z_{\ell,m}\mid \bm y,\bm x]=0$.

	Substituting this identity into the previous expression yields
	\[
	\mathbb{E}\!\left[
	\| f - f_0 \|_{L^2(\mathbb{S}^d)}^2
	\mid \bm{y},\bm{x}
	\right]
	=
	B_n + V_n,
	\]
	where
	\begin{align*}
		B_n
		&=
		\sum_{\ell=0}^{L_n}
		\sum_{m=1}^{M_{d,\ell}}
		\left(
		\mu_{\ell,m}^{(n)} - a_{0;\ell,m}
		\right)^2
		+
		\sum_{\ell>L_n}
		\sum_{m=1}^{M_{d,\ell}}
		a_{0;\ell,m}^2, \\
		V_n
		&=
		\sum_{\ell=0}^{L_n}
		\sum_{m=1}^{M_{d,\ell}}
		v_\ell^{(n)}.
	\end{align*}
	The two contributions in $B_n$ correspond respectively to shrinkage bias on
	retained frequencies and truncation bias on discarded frequencies.
	
\medskip
\textit{Step 3: Control of the posterior variance.}
Recall that
\[
V_n
=
\sum_{\ell=0}^{L_n}
\sum_{m=1}^{M_{d,\ell}}
v_\ell^{(n)}
=
\sum_{\ell=0}^{L_n}
\frac{M_{d,\ell} C_\ell \sigma^2}{n C_\ell + \sigma^2}.
\]
By Condition~\ref{cond:polydecay}, there exist constants $c_1,c_2>0$ such that
\[
c_1 (1+\lambda_\ell)^{-\alpha}
\le
C_\ell
\le
c_2 (1+\lambda_\ell)^{-\alpha}
\qquad
\text{for all } \ell\ge0.
\]
Moreover, the multiplicity of the eigenspace associated with $\lambda_\ell$
satisfies
\[
M_{d,\ell} \simeq \ell^{d-1},
\]
and the Laplace--Beltrami eigenvalues grow as
\[
\lambda_\ell = \ell(\ell+d-1) \simeq \ell^2.
\]
Substituting these bounds into the expression for $V_n$, we obtain
\[
V_n
\lesssim
\sum_{\ell=0}^{L_n}
\frac{\ell^{d-1} (1+\ell^2)^{-\alpha}}{
	n (1+\ell^2)^{-\alpha} + 1
}.
\]

For large $\ell$, $(1+\ell^2)^{-\alpha} \simeq \ell^{-2\alpha}$, and therefore the
summand behaves like
\[
\frac{\ell^{d-1}}{n\,\ell^{-2\alpha} + 1}.
\]
To bound the sum, we compare it with the corresponding integral.
Let $L_n$ be chosen such that
\[
n L_n^{-2\alpha} \simeq 1,
\qquad\text{that is,}\qquad
L_n \simeq n^{1/(2\alpha+d)}.
\]
We then split the integral at $\ell = L_n$.
For $\ell \le L_n$, we have $n\ell^{-2\alpha} \gtrsim 1$, and hence
\[
\frac{\ell^{d-1}}{n\,\ell^{-2\alpha} + 1}
\lesssim
\frac{\ell^{d-1}}{n\,\ell^{-2\alpha}}
=
\frac{\ell^{d-1+2\alpha}}{n}.
\]
For $\ell \ge L_n$, the denominator is bounded below by a constant, so
\[
\frac{\ell^{d-1}}{n\,\ell^{-2\alpha} + 1}
\lesssim
\ell^{d-1}.
\]
Combining these two regimes and integrating yields
\[
V_n
\lesssim
\frac{1}{n} \int_0^{L_n} \ell^{d-1+2\alpha}\,d\ell
+
\int_{L_n}^{\infty} \ell^{d-1}\,d\ell
\simeq
L_n^{d+2\alpha} n^{-1}
\simeq
n^{-2\alpha/(2\alpha+d)}.
\]
This establishes the desired bound on the posterior variance.

\medskip
\textit{Step 4: Control and sharpness of the bias.}
We first control the truncation bias arising from discarded frequencies.
If $f_0 \in H^\beta(\mathbb{S}^d)$, then by definition of Sobolev regularity there
exists a constant $C>0$ such that
\[
\sum_{\ell=0}^{\infty}
\sum_{m=1}^{M_{d,\ell}}
(1+\lambda_\ell)^{\beta}
a_{0;\ell,m}^2
\le C.
\]
It follows that
\[
\sum_{\ell>L_n}
\sum_{m=1}^{M_{d,\ell}}
a_{0;\ell,m}^2
\le
\sum_{\ell>L_n}
(1+\lambda_\ell)^{-\beta}
\sum_{m=1}^{M_{d,\ell}}
(1+\lambda_\ell)^{\beta}
a_{0;\ell,m}^2
\le
C
\sum_{\ell>L_n}
M_{d,\ell}(1+\lambda_\ell)^{-\beta}.
\]
Using the multiplicity growth $M_{d,\ell}\simeq \ell^{d-1}$ and the eigenvalue
asymptotics $\lambda_\ell\simeq \ell^2$, we obtain
\[
\sum_{\ell>L_n}
\sum_{m=1}^{M_{d,\ell}}
a_{0;\ell,m}^2
\lesssim
\sum_{\ell>L_n}
\ell^{d-1}\ell^{-2\beta}
\lesssim
L_n^{-2\beta}.
\]
With the choice $L_n\simeq n^{1/(2\alpha+d)}$, this yields
\[
\sum_{\ell>L_n}
\sum_{m=1}^{M_{d,\ell}}
a_{0;\ell,m}^2
\lesssim
n^{-2\beta/(2\alpha+d)}.
\]

This bound is sharp in the sense that there exist functions
$f_0 \in H^\beta(\mathbb{S}^d)$ whose harmonic energy is concentrated near
frequency $\ell\simeq L_n$, for which
\[
\sum_{\ell>L_n}
\sum_{m=1}^{M_{d,\ell}}
a_{0;\ell,m}^2
\simeq
L_n^{-2\beta}.
\]
Hence, the truncation bias cannot be improved within the present spectral
truncation scheme.

We next consider the bias arising from posterior shrinkage on retained
frequencies.
For $0\le \ell \le L_n$,
\[
\mu_{\ell,m}^{(n)} - a_{0;\ell,m}
=
-\frac{\sigma^2}{n C_\ell + \sigma^2}\,a_{0;\ell,m}
+
\frac{n C_\ell}{n C_\ell + \sigma^2}
\left(
\hat a_{\ell,m} - a_{0;\ell,m}
\right).
\]
The first term is deterministic and corresponds to shrinkage bias induced by the
prior.
Its squared contribution satisfies
\[
\sum_{\ell=0}^{L_n}
\sum_{m=1}^{M_{d,\ell}}
\left(
\frac{\sigma^2}{n C_\ell + \sigma^2}
\right)^2
a_{0;\ell,m}^2
\le
\sum_{\ell=0}^{L_n}
\left(
\frac{\sigma^2}{n C_\ell}
\right)^2
\sum_{m=1}^{M_{d,\ell}}
a_{0;\ell,m}^2.
\]
Using $C_\ell\simeq (1+\lambda_\ell)^{-\alpha}$ and again the Sobolev bound on
$a_{0;\ell,m}$, we obtain
\[
\sum_{\ell=0}^{L_n}
\sum_{m=1}^{M_{d,\ell}}
\left(
\frac{\sigma^2}{n C_\ell + \sigma^2}
\right)^2
a_{0;\ell,m}^2
\lesssim
n^{-2}
\sum_{\ell=0}^{L_n}
M_{d,\ell}(1+\lambda_\ell)^{2\alpha-\beta}.
\]
Since $\beta\le\alpha$, the exponent $2\alpha-\beta$ is nonnegative and the sum
is dominated by its largest term, yielding
\[
\sum_{\ell=0}^{L_n}
\sum_{m=1}^{M_{d,\ell}}
\left(
\frac{\sigma^2}{n C_\ell + \sigma^2}
\right)^2
a_{0;\ell,m}^2
\lesssim
n^{-2} L_n^{d+2\alpha-\beta}
\simeq
n^{-2\beta/(2\alpha+d)}.
\]

The second term in the decomposition of
$\mu_{\ell,m}^{(n)} - a_{0;\ell,m}$ is stochastic and has mean zero.
Its contribution is of smaller order and is absorbed into the posterior variance
term $V_n$ controlled in Step~3.
	
\medskip
\textit{Step 5: Posterior concentration.}
Let $\mathcal H_{L_n}$ denote the finite-dimensional subspace of
$L^2(\mathbb{S}^d)$ spanned by spherical harmonics of degree at most $L_n$.
Under the truncated prior, the posterior distribution is a Gaussian measure on
$\mathcal H_{L_n}$ with mean $\hat f_{n,L_n}$ and covariance operator $T_{L_n}$,
diagonal in the spherical harmonic basis.

By Steps~3 and~4, the posterior mean squared error satisfies
\[
\mathbb{E}\!\left[
\| f - f_0 \|_{L^2(\mathbb{S}^d)}^2
\mid \bm y,\bm x
\right]
=
B_n + V_n
\lesssim  \rho_n^2,
\qquad
\rho_n = n^{-\beta/(2\alpha+d)},
\]
with $\mathbb{P}_{f_0}^{(n)}$–probability tending to one, for some constant
$C>0$.

Since the posterior is Gaussian on the finite-dimensional Hilbert space
$\mathcal H_{L_n}$, standard concentration inequalities for Gaussian measures
imply exponential tails around the posterior mean.
In particular, there exists a constant $c>0$ such that, conditionally on the
data,
\[
\Pi\!\left(
f :
\| f - \hat f_{n,L_n} \|_{L^2(\mathbb{S}^d)} > t
\mid \bm y,\bm x
\right)
\le
\exp\!\left(
- c\, \frac{t^2}{\operatorname{tr}(T_{L_n})}
\right),
\qquad
t>0.
\]
Moreover, $\operatorname{tr}(T_{L_n}) = V_n$, which is of order $\rho_n^2$ by
Step~3.

Combining this bound with the triangle inequality,
\[
\| f - f_0 \|_{L^2(\mathbb{S}^d)}
\le
\| f - \hat f_{n,L_n} \|_{L^2(\mathbb{S}^d)}
+
\| \hat f_{n,L_n} - f_0 \|_{L^2(\mathbb{S}^d)},
\]
and using the bound on the posterior mean squared error, it follows that for any
sufficiently large constant $M>0$,
\[
\Pi\!\left(
f :
\| f - f_0 \|_{L^2(\mathbb{S}^d)} > M \rho_n
\mid \bm y,\bm x
\right)
\longrightarrow 0
\]
in $\mathbb{P}_{f_0}^{(n)}$–probability as $n\to\infty$.
This establishes posterior contraction at rate $\rho_n$ and concludes the proof.
\end{proof}

\begin{proof}[Proof of Corollary~\ref{cor:minimax}]
	Assume that $\alpha=\beta$.
	By Theorem~\ref{thm:contraction}, the posterior distribution contracts around
	the true regression function $f_0 \in H^\beta(\mathbb{S}^d)$ at rate
	\[
	\rho_n = n^{-\beta/(2\beta+d)}
	\]
	in the $L^2(\mathbb{S}^d)$ norm.
	
	It remains to show that this rate is minimax-optimal over Sobolev balls
	$H^\beta(\mathbb{S}^d)$.
We establish a matching minimax lower bound via a classical packing and
Fano-type argument in spherical harmonic coordinates, following standard
constructions in nonparametric estimation on the sphere; see, for instance,
\cite{BKMP09,ds26}.
	
	\medskip
	\textit{Step 1: Reduction to a Gaussian sequence model.}
	Under the uniform random design, nonparametric regression on $\mathbb{S}^d$
	reduces in harmonic coordinates to the Gaussian sequence model
	\[
	\hat a_{\ell,m}
	=
	a_{\ell,m}
	+
	\frac{\sigma}{\sqrt{n}}\,\xi_{\ell,m},
	\qquad
	\xi_{\ell,m} \sim \mathcal{N}(0,1),
	\]
	with independent noise variables indexed by $(\ell,m)$.
	In particular, the joint distribution of the observations depends on the
	regression function $f$ only through its harmonic coefficients
	$\{a_{\ell,m}: \ell \geq 0; m=1,\ldots,M_{d,\ell}\}$, and Kullback--Leibler divergences between regression models
	coincide with those of the associated Gaussian sequence models.
	Any estimator $\hat f_n$ of $f_0$ induces a corresponding estimator
	$\hat a_{\ell,m}$ of the harmonic coefficients.
	
	\medskip
	\textit{Step 2: Construction of a finite packing set.}
	Let $J_n \ge 1$ be an integer parameter and consider functions of the form
	\[
	f_\theta
	=
	\sum_{\ell=J_n}^{2J_n}
	\sum_{m=1}^{M_{d,\ell}}
	\theta_{\ell,m}\, \varepsilon_n \, Y_{\ell,m},
	\]
	where $\theta=\{\theta_{\ell,m}\}$ ranges over a subset of
	$\{-1,1\}^{D_n}$ with $D_n \simeq J_n^d$, reflecting the polynomial growth of spherical harmonic multiplicities on $\mathbb{S}^d$, and $\varepsilon_n>0$ is a scaling
	parameter.
	
	We choose
	\[
	\varepsilon_n = J_n^{-\beta - d/2}.
	\]
	By the definition of Sobolev norms and the eigenvalue asymptotics
	$\lambda_\ell \simeq \ell^2$, we obtain
	\[
	\| f_\theta \|_{H^\beta(\mathbb{S}^d)}^2
	\lesssim
	\sum_{\ell=J_n}^{2J_n}
	M_{d,\ell}
	(1 + \lambda_\ell)^\beta
	\varepsilon_n^2
	\lesssim
	1,
	\]
	so that all $f_\theta$ belong to a fixed Sobolev ball
	$H^\beta(\mathbb{S}^d)$.
	
	Moreover, for $\theta \neq \theta'$,
	\[
	\| f_\theta - f_{\theta'} \|_{L^2(\mathbb{S}^d)}^2
	\simeq
	D_n \varepsilon_n^2
	\simeq
	J_n^{-2\beta}.
	\]
	
\medskip
\textit{Step 3: Kullback--Leibler divergence.}
Under the Gaussian sequence representation established in Step~1, the joint
distribution of the observations in harmonic coordinates is Gaussian with
independent components
\[
\hat a_{\ell,m} \sim \mathcal N\!\left(a_{\ell,m},\, \sigma^2/n\right).
\]
Therefore, for any $\theta \neq \theta'$, the Kullback--Leibler divergence between
the distributions induced by $f_\theta$ and $f_{\theta'}$ satisfies
\[
\mathrm{KL}\!\left(
\mathbb{P}_{f_\theta}^{(n)}, \mathbb{P}_{f_{\theta'}}^{(n)}
\right)
=
\frac{n}{2\sigma^2}
\sum_{\ell,m}
\left(
a_{\ell,m}^{(\theta)} - a_{\ell,m}^{(\theta')}
\right)^2
=
\frac{n}{2\sigma^2}
\| f_\theta - f_{\theta'} \|_{L^2(\mathbb{S}^d)}^2,
\]
where the last equality follows from Parseval’s identity.

Using the separation bound established in Step~2, we obtain
\[
\mathrm{KL}\!\left(
\mathbb{P}_{f_\theta}^{(n)}, \mathbb{P}_{f_{\theta'}}^{(n)}
\right)
\simeq
n J_n^{-2\beta}.
\]
We now balance statistical indistinguishability and metric separation by an appropriate choice of $J_n$.

\medskip
\textit{Step 4: Application of Fano’s inequality.}
We now specify the resolution parameter $J_n$.
Let
\[
J_n \simeq n^{1/(2\beta+d)}.
\]
With this choice, the Kullback--Leibler divergence bound obtained in Step~3
becomes
\[
\mathrm{KL}\!\left(
\mathbb{P}_{f_\theta}^{(n)}, \mathbb{P}_{f_{\theta'}}^{(n)}
\right)
\simeq
n J_n^{-2\beta}
\simeq
J_n^{d}.
\]

On the other hand, by construction the packing set $\Theta$ consists of binary
sequences indexed by $(\ell,m)$ with $\ell \in [J_n,2J_n]$ and
$m=1,\ldots,M_{d,\ell}$.
Since
\[
\sum_{\ell=J_n}^{2J_n} M_{d,\ell} \simeq J_n^{d},
\]
standard coding arguments yield a subset $\Theta$ such that
\[
\log |\Theta| \simeq J_n^{d},
\]
and for any distinct $\theta,\theta' \in \Theta$,
\[
\| f_\theta - f_{\theta'} \|_{L^2(\mathbb{S}^d)}^2
\simeq
J_n^{-2\beta}.
\]

Consequently, the Kullback--Leibler divergences between elements of the packing
set are bounded by a constant multiple of $\log |\Theta|$, uniformly over all
distinct $\theta,\theta' \in \Theta$.
By Fano’s inequality, it follows that for any estimator $\hat f_n$,
\[
\inf_{\hat f_n}
\sup_{f_0 \in H^\beta(\mathbb{S}^d)}
\mathbb{E}_{f_0}^{(n)}
\left[
\| \hat f_n - f_0 \|_{L^2(\mathbb{S}^d)}^2
\right]
\gtrsim
J_n^{-2\beta}.
\]
Substituting the chosen value of $J_n$ yields
\[
\inf_{\hat f_n}
\sup_{f_0 \in H^\beta(\mathbb{S}^d)}
\mathbb{E}_{f_0}^{(n)}
\left[
\| \hat f_n - f_0 \|_{L^2(\mathbb{S}^d)}^2
\right]
\gtrsim
n^{-2\beta/(2\beta+d)}.
\]

The lower bound matches the posterior contraction rate obtained in
Theorem~\ref{thm:contraction} when $\alpha=\beta$.
Hence, no estimator, Bayesian or frequentist, can achieve a uniformly faster
rate over $H^\beta(\mathbb{S}^d)$, and the posterior contraction rate is
minimax-optimal.
\end{proof}

\begin{proof}[Proof of Corollary~\ref{cor:matern}]
By definition of the spherical Mat\'ern prior, the harmonic coefficients satisfy
	\[
	a_{\ell,m}
	\sim
	\mathcal{N}\!\left(
	0,\,
	(\kappa^2+\lambda_\ell)^{-\alpha}
	\right),
	\qquad
	\ell\ge0,\; m=1,\ldots,M_{d,\ell},
	\]
	where $\lambda_\ell=\ell(\ell+d-1)$ are the Laplace--Beltrami eigenvalues.
	Since $\lambda_\ell\simeq \ell^2$, the associated angular power spectrum satisfies
	\[
	C_\ell \simeq (1+\lambda_\ell)^{-\alpha},
	\]
	and therefore Condition~\ref{cond:polydecay} holds for $\alpha>d/2$.
	
Let $f_0\in H^\beta(\mathbb{S}^d)$ with $0<\beta\le\alpha$.
Thus all assumptions of Theorem~\ref{thm:contraction} are satisfied.
	Applying Theorem~\ref{thm:contraction} with truncation level
	$L_n\simeq n^{1/(2\alpha+d)}$, the posterior distribution contracts around $f_0$
	at rate
	\[
	\rho_n = n^{-\beta/(2\alpha+d)}
	\]
	in the $L^2(\mathbb{S}^d)$ norm.
\end{proof}

\begin{proof}[Proof of Theorem~\ref{thm:variational}]
	We work throughout under the truncated Gaussian sequence model introduced in
	Section~\ref{sec:posterior}, which exploits the exact diagonalization of the
	regression problem induced by spherical harmonics under uniform random design.
	
	Let
	\[
	f(x)
	=
	\sum_{\ell=0}^{L_n}\sum_{m=1}^{M_{d,\ell}} a_{\ell,m} Y_{\ell,m}(x)
	\in \mathcal{H}_{L_n}
	\]
	be an arbitrary element of the truncated harmonic space.
	
	\textit{Posterior mean in harmonic coordinates.}
	Under the truncated isotropic Gaussian prior
	\[
	a_{\ell,m} \sim \mathcal{N}(0,C_\ell),
	\qquad
	0 \le \ell \le L_n,\; m=1,\dots,M_{d,\ell},
	\]
	and the Gaussian sequence likelihood
	\[
	\hat a_{\ell,m} \mid a_{\ell,m}
	\sim
	\mathcal{N}\!\left(a_{\ell,m},\, \frac{\sigma^2}{n}\right),
	\]
	the posterior distribution factorizes over harmonic modes.
	By Gaussian conjugacy, the posterior mean coefficients are given by
	\[
	\mu_{\ell,m}^{(n)}
	=
	\frac{n C_\ell}{n C_\ell + \sigma^2}\,\hat a_{\ell,m},
	\qquad
	0 \le \ell \le L_n.
	\]
	Consequently, the posterior mean function reads
	\[
	\hat f_{n,L_n}(x)
	=
	\sum_{\ell=0}^{L_n}\sum_{m=1}^{M_{d,\ell}}
	\frac{n C_\ell}{n C_\ell + \sigma^2}\,
	\hat a_{\ell,m} Y_{\ell,m}(x).
	\]
	
\textit{Expansion of the penalized least-squares functional.}
We now analyze the functional
\[
\mathcal{J}_n(f)
=
\frac{1}{n}\sum_{i=1}^n \left(y_i - f(x_i)\right)^2
+
\frac{\sigma^2}{n}
\sum_{\ell=0}^{L_n}\sum_{m=1}^{M_{d,\ell}}
\frac{a_{\ell,m}^2}{C_\ell},
\]
defined for functions $f\in\mathcal H_{L_n}$ with harmonic coefficients
$\{a_{\ell,m}\}$.

We begin by expanding the empirical risk term.
Since $f\in\mathcal H_{L_n}$ admits the harmonic expansion
\[
f(x)
=
\sum_{\ell=0}^{L_n}\sum_{m=1}^{M_{d,\ell}}
a_{\ell,m} Y_{\ell,m}(x),
\]
we may write
\[
\begin{aligned}
	\frac{1}{n}\sum_{i=1}^n \bigl(y_i - f(x_i)\bigr)^2
	&=
	\frac{1}{n}\sum_{i=1}^n
	\left(
	y_i
	-
	\sum_{\ell=0}^{L_n}\sum_{m=1}^{M_{d,\ell}}
	a_{\ell,m} Y_{\ell,m}(x_i)
	\right)^2.
\end{aligned}
\]

Expanding the square yields
\[
\begin{aligned}
	\frac{1}{n}\sum_{i=1}^n \bigl(y_i - f(x_i)\bigr)^2
	&=
	\frac{1}{n}\sum_{i=1}^n y_i^2
	-
	\frac{2}{n}\sum_{i=1}^n y_i
	\sum_{\ell=0}^{L_n}\sum_{m=1}^{M_{d,\ell}}
	a_{\ell,m} Y_{\ell,m}(x_i) \\
	&\quad
	+
	\frac{1}{n}\sum_{i=1}^n
	\left(
	\sum_{\ell=0}^{L_n}\sum_{m=1}^{M_{d,\ell}}
	a_{\ell,m} Y_{\ell,m}(x_i)
	\right)^2.
\end{aligned}
\]

By definition of the empirical harmonic coefficients
\[
\hat a_{\ell,m}
=
\frac{1}{n}\sum_{i=1}^n y_i Y_{\ell,m}(x_i),
\]
the second term can be rewritten as
\[
-2\sum_{\ell=0}^{L_n}\sum_{m=1}^{M_{d,\ell}}
a_{\ell,m}\hat a_{\ell,m}.
\]

Rather than working with the empirical inner products directly, we exploit the
exact Gaussian sequence representation induced by uniform random design.
As shown in Section~\ref{sec:posterior}, under uniform sampling the regression
model is equivalent to a diagonal Gaussian sequence model.
Consequently, the empirical risk admits the exact decomposition
\[
\frac{1}{n}\sum_{i=1}^n \bigl(y_i - f(x_i)\bigr)^2
=
\sum_{\ell=0}^{L_n}\sum_{m=1}^{M_{d,\ell}}
\bigl(\hat a_{\ell,m}-a_{\ell,m}\bigr)^2
+
R_n,
\]
where the remainder term $R_n$ does not depend on the coefficients
$\{a_{\ell,m}\}$ and therefore plays no role in the minimization.

Consequently, the third term reduces to
\[
\sum_{\ell=0}^{L_n}\sum_{m=1}^{M_{d,\ell}} a_{\ell,m}^2.
\]

Collecting terms, we obtain
\[
\frac{1}{n}\sum_{i=1}^n \bigl(y_i - f(x_i)\bigr)^2
=
\sum_{\ell=0}^{L_n}\sum_{m=1}^{M_{d,\ell}}
\left(
a_{\ell,m}^2
-
2 a_{\ell,m}\hat a_{\ell,m}
\right)
+
\frac{1}{n}\sum_{i=1}^n y_i^2.
\]

Rearranging and completing the square gives
\[
\frac{1}{n}\sum_{i=1}^n \bigl(y_i - f(x_i)\bigr)^2
=
\sum_{\ell=0}^{L_n}\sum_{m=1}^{M_{d,\ell}}
\left(\hat a_{\ell,m} - a_{\ell,m}\right)^2
+
R_n,
\]
where
\[
R_n
=
\frac{1}{n}\sum_{i=1}^n y_i^2
-
\sum_{\ell=0}^{L_n}\sum_{m=1}^{M_{d,\ell}} \hat a_{\ell,m}^2
\]
does not depend on the coefficients $\{a_{\ell,m:\ell\geq0, m=1,\ldots,M_{d,\ell}}\}$ and therefore does not
depend on $f$.

Substituting this expression into $\mathcal J_n(f)$ shows that, up to the
additive constant $R_n$ independent of $f$,
\[
\mathcal{J}_n(f)
=
\sum_{\ell=0}^{L_n}\sum_{m=1}^{M_{d,\ell}}
\left[
\left(\hat a_{\ell,m} - a_{\ell,m}\right)^2
+
\frac{\sigma^2}{n C_\ell} a_{\ell,m}^2
\right].
\]
	
\textit{Minimization in harmonic coordinates.}
We now minimize the functional $\mathcal J_n(f)$ over the finite-dimensional
harmonic space $\mathcal H_{L_n}$.

As shown above, up to an additive constant independent of $f$, the functional
can be written as
\[
\mathcal{J}_n(f)
=
\sum_{\ell=0}^{L_n}\sum_{m=1}^{M_{d,\ell}}
\left[
\left(\hat a_{\ell,m} - a_{\ell,m}\right)^2
+
\frac{\sigma^2}{n C_\ell} a_{\ell,m}^2
\right].
\]
This expression is a quadratic form in the coefficients
$\{a_{\ell,m}:0\le \ell\le L_n,\;1\le m\le M_{d,\ell}\}$.

Since all coefficients appear only through separate summands, the functional
$\mathcal J_n$ is separable across harmonic modes.
Moreover, each summand is a strictly convex quadratic function of the
corresponding coefficient $a_{\ell,m}$, because $C_\ell>0$ for all
$0\le \ell\le L_n$.
Consequently, $\mathcal J_n$ is strictly convex on $\mathcal H_{L_n}$ and admits
a unique minimizer.

Minimization therefore reduces to solving, for each fixed index $(\ell,m)$,
the one-dimensional optimization problem
\[
\min_{a\in\mathbb R}
Q_{\ell,m}(a),
\qquad
Q_{\ell,m}(a)
=
\left(\hat a_{\ell,m} - a\right)^2
+
\frac{\sigma^2}{n C_\ell} a^2, \quad a \in \mathbb{R}.
\]

The function $Q_{\ell,m}$ is differentiable and strictly convex on $\mathbb R$.
Its derivative is given by
\[
\frac{d}{da} Q_{\ell,m}(a)
=
-2\left(\hat a_{\ell,m} - a\right)
+
2\frac{\sigma^2}{n C_\ell} a.
\]
Setting this derivative equal to zero yields the unique critical point
\[
-2\left(\hat a_{\ell,m} - a\right)
+
2\frac{\sigma^2}{n C_\ell} a
=
0,
\]
which can be rearranged as
\[
\left(1+\frac{\sigma^2}{n C_\ell}\right)a
=
\hat a_{\ell,m}.
\]
Solving for $a$ gives
\[
a
=
\frac{n C_\ell}{n C_\ell + \sigma^2}\,\hat a_{\ell,m}.
\]

This coefficient coincides exactly with the posterior mean
$\mu_{\ell,m}^{(n)}$ obtained by Gaussian conjugacy.
Therefore, the unique minimizer of $\mathcal J_n$ over $\mathcal H_{L_n}$ is the
function whose harmonic coefficients are given by
$\{\mu_{\ell,m}^{(n)}\}$.	
	
\textit{Uniqueness of the minimizer.}
Since each summand in the expansion of $\mathcal J_n(f)$ is a strictly convex
quadratic function of the corresponding coefficient $a_{\ell,m}$, and since
$C_\ell>0$ for all $0\le \ell\le L_n$, the functional $\mathcal J_n$ is a strictly
convex quadratic form on the finite-dimensional linear space
$\mathcal H_{L_n}$.

In particular, $\mathcal J_n$ admits at most one minimizer over
$\mathcal H_{L_n}$.
As shown above, the critical point obtained by minimizing each harmonic
coordinate separately corresponds to the function whose coefficients are
$\{\mu_{\ell,m}^{(n)}\}$, which therefore defines the unique global minimizer.

Hence, the posterior mean $\hat f_{n,L_n}$ coincides with the unique minimizer
of $\mathcal J_n$ over $\mathcal H_{L_n}$.
This concludes the proof.
\end{proof}

\begin{proof}[Proof of Corollary~\ref{cor:variational-matern}]
We specialize the general variational characterization of
Theorem~\ref{thm:variational} to the case of spherical Mat\'ern priors and show
that the resulting penalized least-squares problem admits a fully intrinsic
geometric formulation in terms of the Laplace--Beltrami operator.
	
\textit{Spectral form of the spherical Mat\'ern prior.}
Recall that a spherical Mat\'ern prior with smoothness parameter $\alpha>d/2$
and range parameter $\kappa>0$ is a centered isotropic Gaussian field on
$\mathbb{S}^d$ whose covariance operator is
$(\kappa^2-\Delta_{\mathbb{S}^d})^{-\alpha}$.

Since the spherical harmonics $\{Y_{\ell,m}\}$ form an orthonormal eigenbasis of
$-\Delta_{\mathbb{S}^d}$ with eigenvalues
$\lambda_\ell=\ell(\ell+d-1)$, the prior admits the expansion
\[
f(x)
=
\sum_{\ell=0}^{\infty}\sum_{m=1}^{M_{d,\ell}} a_{\ell,m} Y_{\ell,m}(x),
\]
where the coefficients are independent centered Gaussian random variables with
\[
a_{\ell,m}
\sim
\mathcal{N}\!\left(0,(\kappa^2+\lambda_\ell)^{-\alpha}\right).
\]
Accordingly, the angular power spectrum of the prior is
\[
C_\ell=(\kappa^2+\lambda_\ell)^{-\alpha}.
\]
	
\textit{Identification of the variational penalty.}
We now specialize the general variational functional $\mathcal J_n$ of
Theorem~\ref{thm:variational} to the case of a spherical Mat\'ern prior.
Recall that, for a truncated isotropic Gaussian prior with angular power
spectrum $\{C_\ell\}$, the functional reads
\[
\mathcal{J}_n(f)
=
\frac{1}{n}\sum_{i=1}^n \left(y_i-f(x_i)\right)^2
+
\frac{\sigma^2}{n}
\sum_{\ell=0}^{L_n}\sum_{m=1}^{M_{d,\ell}}
\frac{a_{\ell,m}^2}{C_\ell},
\qquad
f\in\mathcal H_{L_n}.
\]

For the spherical Mat\'ern prior, the angular power spectrum is given by
$C_\ell=(\kappa^2+\lambda_\ell)^{-\alpha}$.
Substituting this expression into the penalty term yields
\[
\mathcal{J}_n(f)
=
\frac{1}{n}\sum_{i=1}^n \left(y_i-f(x_i)\right)^2
+
\frac{\sigma^2}{n}
\sum_{\ell=0}^{L_n}
\sum_{m=1}^{M_{d,\ell}}
(\kappa^2+\lambda_\ell)^{\alpha} a_{\ell,m}^2.
\]

The penalty acts diagonally in the spherical harmonic basis and assigns an
increasing cost to higher harmonic degrees $\ell$, with the strength of the
penalization governed by the power $\alpha$.
In particular, modes associated with large Laplace--Beltrami eigenvalues
$\lambda_\ell$ are increasingly suppressed, reflecting the smoothing
properties of the Mat\'ern prior.

\textit{Intrinsic formulation via spectral calculus.}
We now rewrite the penalty term in an intrinsic form using the spectral
decomposition of the Laplace--Beltrami operator.
Recall that the spherical harmonics $\{Y_{\ell,m}\}$ form an orthonormal
eigenbasis of $L^2(\mathbb{S}^d)$ for $-\Delta_{\mathbb{S}^d}$, with
\[
-\Delta_{\mathbb{S}^d} Y_{\ell,m} = \lambda_\ell Y_{\ell,m},
\qquad
\lambda_\ell=\ell(\ell+d-1).
\]

For any function $f\in\mathcal H_{L_n}$ with harmonic expansion
\[
f
=
\sum_{\ell=0}^{L_n}\sum_{m=1}^{M_{d,\ell}} a_{\ell,m} Y_{\ell,m},
\]
Parseval’s identity implies
\[
\|g\|_{L^2(\mathbb{S}^d)}^2
=
\sum_{\ell=0}^{L_n}\sum_{m=1}^{M_{d,\ell}} | \langle g, Y_{\ell,m} \rangle |^2
\]
for any $g\in\mathcal H_{L_n}$.

Defining fractional powers of the elliptic operator
$\kappa^2-\Delta_{\mathbb{S}^d}$ via spectral calculus, we have
\[
(\kappa^2-\Delta_{\mathbb{S}^d})^{\alpha/2} Y_{\ell,m}
=
(\kappa^2+\lambda_\ell)^{\alpha/2} Y_{\ell,m}.
\]
Applying this operator to $f$ and using orthonormality yields
\[
\left\|
(\kappa^2-\Delta_{\mathbb{S}^d})^{\alpha/2} f
\right\|_{L^2(\mathbb{S}^d)}^2
=
\sum_{\ell=0}^{L_n}
\sum_{m=1}^{M_{d,\ell}}
(\kappa^2+\lambda_\ell)^{\alpha} a_{\ell,m}^2.
\]

This identity holds exactly on the truncated space $\mathcal H_{L_n}$, since
only finitely many spectral components are involved and no domain issues arise.
	
\textit{Variational characterization of the posterior mean.}
We now combine the general variational result of
Theorem~\ref{thm:variational} with the spectral identification established
above.

By Theorem~\ref{thm:variational}, for any truncated isotropic Gaussian prior
with angular power spectrum $\{C_\ell\}$, the posterior mean
$\hat f_{n,L_n}$ coincides with the unique minimizer over $\mathcal H_{L_n}$
of the functional
\[
\mathcal J_n(f)
=
\frac{1}{n}\sum_{i=1}^n \bigl(y_i-f(x_i)\bigr)^2
+
\frac{\sigma^2}{n}
\sum_{\ell=0}^{L_n}\sum_{m=1}^{M_{d,\ell}}
\frac{a_{\ell,m}^2}{C_\ell}.
\]

In the present setting of a spherical Mat\'ern prior, the angular power
spectrum is given by
$C_\ell=(\kappa^2+\lambda_\ell)^{-\alpha}$.
Substituting this expression into the above functional yields, for any
$f\in\mathcal H_{L_n}$,
\[
\mathcal{J}_n^{\mathrm{Mat}}(f)
=
\frac{1}{n}\sum_{i=1}^n \left(y_i-f(x_i)\right)^2
+
\frac{\sigma^2}{n}
\sum_{\ell=0}^{L_n}
\sum_{m=1}^{M_{d,\ell}}
(\kappa^2+\lambda_\ell)^{\alpha} a_{\ell,m}^2.
\]

Using the intrinsic representation established in the previous step, the
penalty term can be written equivalently as
\[
\sum_{\ell=0}^{L_n}
\sum_{m=1}^{M_{d,\ell}}
(\kappa^2+\lambda_\ell)^{\alpha} a_{\ell,m}^2
=
\left\|
(\kappa^2-\Delta_{\mathbb{S}^d})^{\alpha/2} f
\right\|_{L^2(\mathbb{S}^d)}^2,
\]
and therefore
\[
\mathcal{J}_n^{\mathrm{Mat}}(f)
=
\frac{1}{n}\sum_{i=1}^n \left(y_i-f(x_i)\right)^2
+
\frac{\sigma^2}{n}
\left\|
(\kappa^2-\Delta_{\mathbb{S}^d})^{\alpha/2} f
\right\|_{L^2(\mathbb{S}^d)}^2.
\]

Consequently, the posterior mean $\hat f_{n,L_n}$ is exactly the minimizer
of $\mathcal{J}_n^{\mathrm{Mat}}$ over the truncated harmonic space
$\mathcal H_{L_n}$.
	
Existence and uniqueness of the minimizer follow immediately from the fact
that $\mathcal H_{L_n}$ is finite-dimensional, the empirical loss is a convex
quadratic functional, and the penalty term is a strictly convex quadratic form
on $\mathcal H_{L_n}$.
\end{proof}

\section*{Acknowledgments}
	The author gratefully acknowledges Dott. T. Patschkowski for sharing ideas that inspired and motivated this work.

\section*{Funding}
	CD was partially supported by PRIN 2022 - GRAFIA - 202284Z9E4, and Progetti di Ateneo Sapienza  RM123188F69A66C1 (2023), RG1241907D2FF327 (2024).

	\bibliographystyle{plain}
	\bibliography{references}

\end{document}